\documentclass[letterpaper,11pt,leqno,oneside]{article}
\usepackage{latexsym,amssymb,amscd,hyperref,graphicx,subfigure,, epsfig, amsmath, amsfonts,subfigure,graphicx,float,titlesec,amsthm,color,sectsty,mathrsfs}

\usepackage[left=1.75cm, top=2.7cm, bottom=2.7cm,right=1.75cm]{geometry}

\def\newblock{\hskip .11em plus .33em minus .07em}
\numberwithin{equation}{section}

\newtheorem{thrm}{Theorem}[section]
\newtheorem{coro}[thrm]{Corollary}
\newtheorem{pro}[thrm]{Proposition}
\newtheorem{rem}[thrm]{Remark}
\newtheorem{lem}[thrm]{Lemma}

\def\mbR{\mathbb{R}}
\def\mbP{\mathbb{P}}
\def\mbE{\mathbb{E}}
\def\mbR{\mathbb{R}}

\renewcommand\thefigure{\arabic{section}.\arabic{figure}}
\renewcommand{\thesubfigure}{\thefigure(\alph{subfigure})}
\makeatletter
\renewcommand{\p@subfigure}{}
\renewcommand{\@thesubfigure}{\thesubfigure:\hskip\subfiglabelskip}
\makeatother

\definecolor{DRed}{rgb}{0.4, 0.0, 0.0}

\definecolor{Green}{rgb}{0.0, 0.6, 0.0}

\definecolor{Black}{rgb}{0, 0, 0}

\definecolor{Orange}{rgb}{1, 0.4, 0}

\definecolor{Brown}{rgb}{0.59, 0.29, 0}

\definecolor{Purple}{rgb}{1, 0, 1}

\title{Small-time expansions for state-dependent local jump-diffusion models with infinite jump activity}

\author{Jos\'e E. Figueroa-L\'{o}pez\thanks{Department of Mathematics, Washington University in St. Louis, St. Louis, MO 63130,  USA ({\tt figueroa@math.wustl.edu}). {Research supported in part by the {NSF Grants: DMS-1149692.}}}
\and Yankeng Luo\thanks{Department of Mathematics and Applied Mathematics, Virginia Commonwealth University,  Richmond, VA 23284,  USA ({\tt yluo@vcu.edu}).}}

\date{\today}

\begin{document}
\maketitle

\begin{abstract}
In this article, we consider a Markov process $\{X_{t}\}_{t\geqslant0}$, starting from $x\in\mbR$ and solving a stochastic differential equation, which is driven by a Brownian motion and an independent pure jump component exhibiting state-dependent jump intensity and infinite jump activity. A second order expansion is derived for the tail probability $\mbP[X_t\geqslant x+y]$ in small time $t$, for $y>0$. As an application of this expansion and a suitable change of the underlying probability measure, a second order expansion, near expiration, for out-of-the-money European call option prices is obtained when the underlying stock price is modeled as the exponential of the jump-diffusion process $\{X_{t}\}_{t\geqslant0}$ under the risk-neutral probability measure.

\vspace{0.3 cm}
\noindent\textbf{Keywords and phrases}: Short-time asymptotics; local jump-diffusion Markov models; stochastic differential equations with jumps; option pricing.
\end{abstract}

\section{Introduction}
In this work we consider a Markov process $X:=\{X_{t}\}_{t\geqslant0}$ with infinitesimal generator of the form
\begin{equation}\label{InfGenX}
	L f(x)=b(x)f'(x)
+\frac{\sigma^2(x)}{2}f''(x)
+\int_{\mbR_{0}}\left(f(x+\gamma(x,r))-f(x)-\mathbf{1}_{\{|r|\leqslant1\}}\gamma(x,r)f'(x)\right)\nu(x,r)dr,
\end{equation}
where $\mbR_{0}:=\mbR\backslash\{0\}$ and $b$, $\sigma$, $\gamma$, and $\nu$ are deterministic function satisfying appropriate conditions for the existence of such a process (see below for further details). Broadly, $X$ can be defined in terms of a stochastic differential equation (SDE) of the form:
\[
	dX_{t}=b(X_{t})dt+\sigma(X_{t})dW_{t}+d J_{t},
\]
where $W:=\{W_{t}\}_{t\geq{}0}$ is a Wiener process and $J:=\{J_{t}\}_{t\geq{}0}$ is an independent pure-jump process, whose jump behavior is dictated by $\nu$ and $\gamma$ as follows:
\begin{align}\nonumber
   \mbE\big[\#\{s\in[t,t+\delta]:\Delta X_s\in(a,b)\}\big]
   &=\mbE\big[\#\{s\in[t,t+\delta]:\Delta J_s\in(a,b)\}\big]\\
   &=\mbE\left[\int_t^{t+\delta}\int {\bf 1}_{\{\gamma(X_{s^-},r)\in (a,b)\}}\nu(X_{s^-},r)drds\right],\label{InstyJmps0}
 \end{align}
 for any $t\in(0,\infty)$, $\delta>0$, and $(a,b)\in\mbR\backslash\{0\}$. Intuitively, (\ref{InstyJmps0}) tells us that the jump intensity of the process ``near" time $t$ depends on its state immediately before $t$ via the function $\nu$ in that if  $\nu(X_{t^-},r)$ is large (small), then we expect a higher (lower) intensity of jumps immediately after time $t$. In the particular case of $\gamma(x,r)\equiv{}r$, (\ref{InstyJmps0}) reduces to
\begin{equation}\label{PrtCse0}
  \mbE\big[\#\{s\in[t,t+\delta]:\Delta X_s\in(a,b)\}\big]=\mbE\left[\int_t^{t+\delta}\int_{a}^{b}\nu(X_{s^-},r)drds\right],
\end{equation}
and $\nu\left(X_{t^-},r\right)$ has the usual interpretation of a stochastic jump intensity as defined in, e.g., \cite{Bremaud} and \cite{Daley}. That is, $\nu(x,r)$ measures the expected number of jumps, per unit time, with size near $r$ when the process is at state $x$. State-dependent jump behavior as described above is an important feature that offers greater modeling flexibility to other commonly studied jump processes. For several applications we refer the reader to \cite{Johnes}, \cite{Kou}, \cite{Merener}, \cite{Glasserman}, \cite{Duffie},  \cite{Pan}, \cite{Lorig},  and \cite{Yu}. 

The generator (\ref{InfGenX}) covers a wide range of processes. For a L\'evy processes, $b$ and $\sigma$ are constants, $\gamma(x,r)=r$, and $\nu(x,r)=h(r)$, for a L\'evy density $h:\mbR\backslash\{0\}\to[0,\infty)$ (i.e., $\int(x^{2}\wedge 1)h(x)dx<\infty$). When {we simply have} $\nu(x,r)=h(r)$, we recover the class of (local) jump-diffusion models studied in \cite{Cheng}. In that case, $X$ can be constructed as
\begin{equation}\label{SDEIntroduction}
	X_{t}=x+\int_{0}^{t}{b}(X_{s})ds+\int_{0}^{t}\sigma(X_{s})dW_{s}+\sum_{s\in(0,t]:|\Delta Z_{u}|\geq{}1}\gamma(X_{s^{-}},\Delta Z_{s})+\sum^{c}_{s\in(0,t]:0<|\Delta Z_{u}|\leq{}1}\gamma(X_{s^{-}},\Delta Z_{s}), 
\end{equation}
where $Z$ is a L\'evy process with L\'evy density $h$ and $\displaystyle{\sum^{c}}$ denotes the compensated Poisson sum of the terms therein. The case of $\nu(x,r)=\lambda(x)p(r)$ with $\int p(r)dr=1$ has been studied in \cite{Yu}. More generally, if $\lambda(x):=\int \nu(x,r)dr<\infty$ is locally bounded, we can construct the process $X$ as 
\[
	X_{t}=x+\int_{0}^{t}\bar{b}(X_{s})ds+\int_{0}^{t}\sigma(X_{s})dW_{s}+\sum_{\tau_{i}\leq{}t}\gamma(X_{\tau_{i}^{-}},\xi_{i}), 
\]
where $\bar{b}(x)=b(x)-\int\mathbf{1}_{\{|r|\leqslant1\}}\gamma(x,r)\nu(x,r)dr$, $0<\tau_{1}<\tau_{2}<\dots$ is a point process on $\mbR_{+}$ with stochastic intensity $\{\lambda(X_{s^{-}})\}_{s\geq{}0}$ and, conditionally on $X_{\tau_{i}^{-}}=z$, $\xi_{i}$ has density $\nu(z,\cdot)/\lambda(z)$, independently of any other process. 

Unlike the just described processes with finite jump activity (i.e., finitely many jumps during any bounded time period), infinite jump activity ({IJA}) processes are important not only from a mathematical but also practical point of view. This is especially true in financial applications where several statistical tests, based on high-frequency observations, have supported the latter feature. In this work, we study a class of IJA processes (i.e., $\int \nu(x,r)dr=\infty$) that arises as the thinning of a local jump-diffusion process driven by a L\'evy process with stable-like small jump behavior.

An important obstacle for the application of the class of processes described above arise from the lack of tractable transition distributions. In this regard, one stream of the literature has focused on numerical methods for the computational simulation of the process (see, e.g., \cite{Glasserman}  and \cite{Mordecki} and references therein), which in turn can be used to estimate different distributional features via Monte Carlo methods. Another stream of the literature has developed approximations for those distributional properties under different asymptotic regimes. Approximations in short-time are of particular relevance due to their wide range of applications such as statistical estimation and simulation methods. The latter approach was further developed by Figueroa-L\'opez, et al.  \cite{Cheng}, where a second order expansion, in small time, for the tail probability of the process (\ref{SDEIntroduction}) was developed. Unfortunately, there is almost no work dealing with short-time approximation methods for the distributional properties of state-dependent jump-diffusion. An exception is \cite{Yu}, where a short-time approximation scheme for the transition densities was developed in the special case that $\nu(x,r)=\lambda(x)p(r)$ with $\int p(r)dr=1$.

In this article, we generalize the result in \cite{Cheng} by developing a second-order expansion, in small time $t$, for the tail probability $\mathbb{P}\left(X_{t}\geq{}x+y\right)$  {($y>0$) of a state-dependent jump-diffusion process $X$ with generator (\ref{InfGenX}) and initial value $x\in\mbR$}. One of the main motivations for considering both the infinitesimal generator (\ref{InfGenX}) and the tail probabilities is their role in the evaluation of out-of-the-money (OTM) European call option prices when the underlying stock price is modeled as the exponential of a state-dependent jump-diffusion $X$ under a risk-neutral probability measure:
\begin{equation}\label{OTMOP0}
	\Pi_{t}=\mbE\left[\left(e^{X_{t}}-e^{\kappa}\right)_{+}\right],
\end{equation}
where $X$ has initial value $0$ and the log-moneyness $\kappa$ is such that $\kappa>0$. An appealing method for evaluating (\ref{OTMOP0}) is to consider the so-called share measure $\mbP^{\#}$ (cf. \cite{CM09}), defined as $d\mbP^{\#}=e^{X_{t}}d\mbP$. In that case, the following neat representation for $\Pi_{t}$ in terms of two tail probabilities holds:
\begin{equation}\label{OTMOPDec}
	\Pi_{t}=\mbE\left[\left(e^{X_{t}}-e^{\kappa}\right){\bf 1}_{X_{t}\geq{}\kappa}\right]=\mbE\left[e^{X_{t}}{\bf 1}_{X_{t}\geq{}\kappa}\right]-e^{\kappa}\mbP\left[X_{t}\geq{}\kappa\right]=\mbP^{\#}\left[X_{t}\geq{}\kappa\right]-e^{\kappa}\mbP\left[X_{t}\geq{}\kappa\right].
\end{equation}
When $X$ is a L\'evy process, it turns out that the law of $X$ under $\mbP^{\#}$ is again L\'evy, albeit with different L\'evy triplet, and the same small-time expansion for the tail probability can be used to deal with the two terms in (\ref{OTMOPDec}) and obtain an expansion for option price $\Pi_{t}$ (cf. \cite{Martin}). For a general process $X$ with generator (\ref{InfGenX}), the law of $X$ under $\mbP^{\#}$ is still Markovian with the same generator form (\ref{InfGenX}) but replacing $\nu$ and $b$ with 
\begin{align}\label{bbound20}
\nu^{\#}(x,r):=e^{\gamma(x,r)}\nu(x,r),\quad 
b^{\#}(x)&=\frac{\sigma^2(x)}{2}
-\int\left(e^{\gamma(x,r)}-1-e^{\gamma(x,r)}\mathbf{1}_{|r|\leqslant1}\gamma(x,r)\right)\nu(x,r)dr,
\end{align}
respectively. Note that, in particular, even if the jump dynamics of $X$ were not state-dependent (i.e., $\nu(x,r)=\nu(r)$),  $X$ would still be state-dependent under $\mbP^{\#}$ (i.e., $\nu^{\#}$ would depend on $r$). This motivates us to study at once the unifying framework (\ref{InfGenX}), which, as explained above, is important in their own right.

Our main results explicitly quantify the effect of the jump state-dependence in the leading and second order terms of the tail probability and OTM call option premium in a small-time setting. Hence, for instance, when $\gamma(x,r)=r$ (in which case $\nu(x,r)$ really measures the intensity of jumps with size near $r$ when the process' state is $x$), the effect of a positive constant drift $b$ is to change the probability of a large move of more than $y$ within time $t$ by
\begin{equation}\label{CD0}
	\frac{t^{2}}{2}b\left(2\nu(x,y)+\int_{y}^{\infty}\frac{\partial\nu(x,r)}{\partial x}dr\right)\left(1+o(1)\right), \quad (t\to{}0).
\end{equation}
Similarly,
a nonzero constant volatility $\sigma$ will change the probability of a large move of more than $y$ by
\begin{equation}\label{Csigma1}
 \frac{t^{2}}{2}\sigma^{2}\left(-\frac{\partial \nu(x,y)}{\partial y}+\frac{1}{2}\int_{y}^{\infty}\frac{\partial^{2}\nu(x,r)}{\partial x^{2}}dr\right)\left(1+o(1)\right), \quad (t\to{}0),
\end{equation}
and the OTM call option premium (\ref{OTMOP0}) by
\begin{equation}\label{Csigma2}	
{\frac{t^2\sigma^2}{2}\bigg[e^k\nu(0,k)+\int_{k}^{\infty}\frac{e^r+e^k}{2}\left.\frac{\partial \nu(x,r)}{\partial x}\right|_{x=0}dr+\int_{k}^{\infty}\frac{e^r-e^k}{2}\left.\frac{\partial^{2}\nu(x,r)}{\partial x^{2}}\right|_{x=0}dr\bigg](1+o(1))}, \quad (t\to0),
\end{equation}
which extends a result in \cite{Martin}  for exponential L\'evy models.

Our approach to obtain the expansion of the tail probabilities of $X$ follows along the same lines of \cite{Cheng} building on a small-large jump decomposition of the process $X$, similar to that introduced by \cite{Leandre}. {An essential ingredient of this method is to show that the small-jump component is a diffeomorphism, which presents some interesting subtleties in} the considered model. More concretely, to {obtain} the latter property under the new state-dependent jump structure, we prove that the state-dependent small jump component has the same law as a regular (i.e., state-free) jump-diffusion with sufficiently {regular} coefficients (see Section \ref{weaksol} for details). Another novelty in this work is a new, much simpler, proof for estimating the tail probability of the so-called one-jump process (the process conditionally on having one big jump at a specified time) based on the just mentioned diffeomorphism {property, time-reversibility, and a suitable application of an iterated Dynkin's formula} (see Lemma \ref{counterpartA.1} for details).

The paper is organized as follows. In Section \ref{setup}, we formally define the model together with the standing assumptions throughout the article. Section \ref{notation&definition} introduces some needed notations and probabilistic tools. The key ingredients to obtain the main result are presented in Sections \ref{weaksol} and \ref{importantlemma}. The second-order expansion for the tail probability is presented in Section \ref{secExp}, while the application thereof to option pricing is developed in Section \ref{option}. Finally, Section \ref{example} presents a numerical example to illustrate the performance of the expansions. To this end, we also develop a simulation method for our state-dependent local jumps-diffusion model based on a suitable diffusive approximation of the small jump component of the process. Finally, all the proofs are deferred to several Appendix sections.

\textbf{General Notation:} Throughout, given an Euclidean domain $E$, $\bar{C}^{k}_{b}(E)$ (resp., $C^{k}_{b}(E)$)
represents the class of $k^{th}$-times differentiable (resp., bounded and $k^{th}$-times differentiable) functions $f:E\to\mbR$ with continuous and bounded partial derivatives of order $n=1,\dots, k$. In particular, $ C^{k}_{b}(E)\subset \bar{C}^{k}_{b}(E)$. Also, $\partial_{i}$ and $\partial_i^n$ respectively denote the derivative and the n-th order partial derivative operator with respect to the i-th variable of a multivariate function.

\section{Setup and assumptions}\label{setup}

In this section, we give a construction of the process of interest and {establish} the assumptions needed throughout. As mentioned in the introduction, we want to consider an infinite-jump activity Markov process $X$ with infinitesimal generator (\ref{InfGenX}). We use a thinning technique for the construction of the jump structure of $X$ based on the jumps of a suitable L\'evy process $Z$. To this end, we imposed the following assumption: 
\renewcommand{\labelenumi}{(\theenumi)}
\renewcommand{\theenumi}{S\arabic{enumi}}
\renewcommand{\labelenumii}{(\theenumii)}
\renewcommand{\theenumii}{\roman{enumii}}
\begin{enumerate}
\item\label{conditionnu}
    \renewcommand{\theenumi}{}
    \begin{enumerate}
      \item\label{hExistence} There exists a L\'evy density $h$ dominating the jump intensity function $\nu:\mbR\times\mbR_0\to(0,\infty)$; i.e., $\nu(x,r)\leqslant{}h(r)$, for all $(x,r)\in\mbR\times\mbR_0$; 
      \item\label{Smothnuandbarnu} We also assume that  
      $\bar{\nu}(x,r):=\nu(x,r)/h(r)$ and $h$ are ${C}^4_{b}(\mbR\times[\epsilon,\epsilon]^{c})$ and ${C}^4_{b}([\epsilon,\epsilon]^{c})$, respectively, for any $\epsilon>0$ and, furthermore, 
\[
	\liminf_{r\to{}0^{\pm}}\inf_{x\in\mbR}\bar\nu(x,r)>0,\quad 
		\limsup_{r\to{}0^{\pm}}\sup_{x\in\mbR}|r \partial_{2}\bar\nu(x,r)|<\infty, 
		\quad 
		\limsup_{r\to{}0^{\pm}}\sup_{x\in\mbR}|\partial_{1}^{i}\bar\nu(x,r)|<\infty,\quad i=0,1,2.
\]
      \end{enumerate}
\end{enumerate}
\renewcommand{\labelenumi}{\theenumi.}
\renewcommand{\theenumi}{\arabic{enumi}}

A necessary and sufficient condition for the first requirement in the previous condition is that $\bar{h}(r):=\sup_{x}\nu(x,r)$ is a L\'evy density. In that case, we can take $h=\bar{h}$; however, in applications it may be more convenient to choose another $h$ whose associate L\'evy process can be simulated more easily. The second condition therein imposes some regularity requirements. Note that indeed the derivatives of $\nu$ do appear in the expansions of the tail probability and OTM call option premium (see, e.g., (\ref{CD0})-(\ref{Csigma2})), which lead us to believe that some smoothness properties on $\nu$ {are} needed.

 We are now ready to give the construction for $X$. Throughout, we consider a filtered probability space $\left(\Omega,\mathcal{F}, \mathbb{F}:=\{\mathcal{F}_t\}_{t\geqslant0},\mbP\right)$ equipped with a standard Brownian motion $\{W_t\}_{t\geqslant0}$ and an independent Poisson random measure $p(dt,dr,du)$ on $\mbR_+\times E:=\mbR_+\times\mbR_0\times(0,1)$ with mean measure  $dt\, h(r)dr\,du$. The compensated Poisson measure of $p$ is denoted by $\bar{p}(ds,dr,du):=p(ds,dr,du)-ds\, h(r)dr\,du$. Then, under the condition (\ref{conditionnu}-i), we have the following construction for the process $X^{(x)}:=\{X^{(x)}_t\}_{t\geqslant0}$:
\begin{equation}\label{SDEpq}
\begin{aligned}
  X_t^{(x)}&=x+\int_0^t{}b(X_s^{(x)})ds
  +\int_0^t\sigma(X_s^{(x)})dW_s
  +\int_0^t\int_{E}\mathbf{1}_{|r|>1}\gamma(X_{s^-}^{(x)},r)\theta(X_{s^-}^{(x)},r,u)p(ds,dr,du)\\
  &\quad+\int_0^t\int_{E}\mathbf{1}_{|r|\leqslant1}\gamma(X_{s^-}^{(x)},r)\theta(X_{s^-}^{(x)},r,u)\bar{p}(ds,dr,du),
\end{aligned}
\end{equation}
where $\theta:\mbR\times\mbR_0\times(0,1)$  is a thinning function that takes the form
\begin{equation}\label{DfnThetafrst}
 	\theta(x,r,u):=\mathbf{1}_{\left\{u<\frac{\nu(x,r)}{h(r)}\right\}}.
\end{equation}
Upon the existence and uniqueness of the solution to (\ref{SDEpq}), the solution process $X^{(x)}$ would be Markovian with infinitesimal generator (\ref{InfGenX}). The following conditions on $b$, $\sigma$, and $\gamma$ guarantee the well-posedness of (\ref{SDEpq}) (as proved in Lemma \ref{Y0diffeomorphism} below) and other needed features of the process:
\renewcommand{\labelenumi}{(\theenumi)}
\renewcommand{\theenumi}{S\arabic{enumi}}
\renewcommand{\labelenumii}{(\theenumii)}
\renewcommand{\theenumii}{\roman{enumii}}
\begin{enumerate}
\setcounter{enumi}{1}
\item\label{Cinfinity}
    \renewcommand{\theenumi}{}
    \begin{enumerate}
      \item\label{bsigma} The functions $b:\mbR\to\mbR$ and $\sigma:\mbR\to\mbR$ belong to ${C}^{4}_{b}(\mathbb{R})$;
      \item\label{sigma>eta} There exists a constant $\eta>0$ such that $\sigma(x)\geqslant\eta$ for all $x\in\mbR$;
    \end{enumerate}
\renewcommand{\theenumi}{S\arabic{enumi}}
\item\label{conditiongamma}The function $\gamma(x,r):\mbR\times\mbR\to\mbR$ satisfies the following conditions:
    \renewcommand{\theenumi}{}
    \begin{enumerate}
        \item\label{gamma0=0} It belongs to $\bar{C}^{4}_{b}(\mbR\times\mbR)$, and $\gamma(x,0)=0$ for all $x\in\mbR$;
        \item\label{partial2gamma>eta} For all $x,r\in\mbR$, $|\partial_2\gamma(x,r)|>\eta$  for some constant $\eta>0$;
        \item\label{1+partial1gamma} For all $x,r\in\mbR$, $|1+\partial_1\gamma(x,r)|>\eta$ for some constant $\eta>0$;
    \end{enumerate}
\end{enumerate}
\renewcommand{\labelenumi}{\theenumi.}
\renewcommand{\theenumi}{\arabic{enumi}}

\begin{rem}\label{remarkonS}Some remarks are in order regarding the above conditions:
\begin{enumerate}
\item The condition (\ref{Cinfinity}-\ref{sigma>eta}) is a standard non-degeneracy condition that is also imposed {in \cite{Cheng}}.

\item The condition (\ref{conditiongamma}-\ref{gamma0=0}) implies that, for each $\varepsilon>0$ and $i=0,\dots,3$, $|\partial_1^i\gamma(x,r)|\leqslant{}C_{\varepsilon}|r|$ for any $|r|\leqslant\varepsilon$ and some constant $C_{\varepsilon}<\infty$. 

\item The condition (\ref{conditiongamma}-\ref{partial2gamma>eta}) implies that the mapping $r\mapsto\gamma(x,r)$ is strictly monotone with range $(-\infty,\infty)$ and, hence, admits an inverse, denoted by $\gamma^{-1}(x,r)$ hereafter. Without loss of generality, throughout we assume that $r\mapsto\gamma(x,r)$ is \emph{strictly monotone increasing}. 

\item Condition (\ref{conditiongamma}-\ref{1+partial1gamma}) is essential for the mapping $x\to X_{t}^{(x)}$ to be a {diffeomorphism (see Lemma \ref{Y0diffeomorphism} below)}.

\end{enumerate}
\end{rem}
It is worth pointing out that Figueroa-Lopez et al.~\cite{Cheng} imposed stronger regularity to the coefficients of their SDE. However, we shall see that most results therein are still valid under {the milder regularities imposed in the} present manuscript.

\begin{rem}\label{RemEquivTM}
As mentioned in the introduction, in the case that $\nu(x,r)= h(r)$, for a L\'evy density $h$, we recover the model (\ref{SDEIntroduction}), which was studied in \cite{Cheng}. Even though it is not evident, under the condition (\ref{conditiongamma}-\ref{partial2gamma>eta}), the law of the process (\ref{SDEpq}) is actually equivalent to that of a process of the form (\ref{SDEIntroduction}) with suitably chosen coefficient functions $\gamma$ and $b$  (see Remark \ref{EILTP} below for more details). However, the resulting {$\gamma$} is relatively intractable and {does not meet the regularity} conditions of \cite{Cheng} for the second-order expansion therein to be applied directly. Furthermore, there are two other important reasons for directly considering the process (\ref{SDEpq}). First, the process (\ref{SDEpq}) allows the direct modeling and clearer interpretation of the intensity of jumps via the parameter $\nu(x,r)$, which is somehow hidden inside the function $\gamma$  in the model (\ref{SDEIntroduction}). Second, as already mentioned in the introduction, in order to develop the small-time expansion of out-of-the-money option prices, one needs to deal with processes having the most general generator (\ref{InfGenX}), even if the original process $X$ is of the form (\ref{SDEIntroduction}).
\end{rem}

Our final assumption is probably the less intuitive. As mentioned in the introduction, there are two key ingredients in our approach to obtain the small-time expansion of the tail probability of $X$. First, we use a small-large jump decomposition of the process $X$. Second, we need that the small-jump component is a diffeomorphism. To this end, we use the equivalence of the resulting state-dependent small-jump model to a state-free jump diffusion process of the form (\ref{SDEIntroduction}), which is possible in light of the above Condition (\ref{conditiongamma}-\ref{partial2gamma>eta}), as described in the previous Remark \ref{RemEquivTM}. {However, to} conclude that a model of the form (\ref{SDEIntroduction}) is a diffeomorphism, we need some {regularity conditions on its coefficients. The main goal of the the following relatively mild condition is to establish such conditions} (see Proposition \ref{deltaproperty} below). We refer to the Remark \ref{WhyS4} below for further discussion and possible relaxation of this condition.
\renewcommand{\labelenumi}{(\theenumi)}
\renewcommand{\theenumi}{S\arabic{enumi}}
\renewcommand{\labelenumii}{(\theenumii)}
\renewcommand{\theenumii}{\roman{enumii}}
\begin{enumerate}
\setcounter{enumi}{3}
\renewcommand{\theenumi}{S\arabic{enumi}}
\item\label{conditionh} The L\'{e}vy density h introduced in the Condition (\ref{conditionnu}) is such that, for some $\alpha\in(0,2)$,  $g(r):=h(r)|r|^{\alpha+1}$ is differentiable in {$(-\epsilon_{0},0)\cup(0,\epsilon_{0})$, for some $\epsilon_{0}>0$}, and 
\[
	{\liminf_{r\to{}0^{\pm}}g(r)>0,\qquad 
	\limsup_{r\to{}0^{\pm}}g(r)<\infty,\qquad
	\limsup_{r\to{}0^{\pm}}|rg'(r)|<\infty.}
\]
\end{enumerate}
\renewcommand{\labelenumi}{\theenumi.}
\renewcommand{\theenumi}{\arabic{enumi}}
`
\section{Some needed notations and preliminary results}\label{notation&definition}

Let $Z:=\{Z_t\}_{t\geqslant0}$ be a pure-jump L\'{e}vy process with L\'{e}vy triplet $(0,h(r)\,dr,0)$, for the truncation function $\mathbf{1}_{|r|\leqslant1}$. The jump measure of the process $Z$ is denoted by $q(dt,dr):=\#\{(t,\Delta{}Z_t)\in{}dt\times{}dr:\Delta{}Z_t\neq0\}=\sum\delta_{(T_i,R_i)}(dt\times{}dr)$, where $\{(T_i,R_i)\}$ are the atoms of the measure $q(dt,dr)$. In that case, the Poisson random measure  $p(dt,dr,du)$ in (\ref{SDEpq}) has the same distribution  as a marked point process on $\{(T_i,R_i)\}$, with marks $\{U_i\}_{i\geqslant1}$ being a random sample from a standard uniform distribution on $(0,1)$.

In the sequel, the process $X^{(x)}$ is decomposed, in law, into a process with small jumps and an independent process of finite jump activity. To formally define the small-jump component of $X$, we first need to introduce a suitable construction for the Poisson random measure $p$.
For any $\varepsilon\in(0,1)$, let 
\[
	h_{\varepsilon}(r):=\phi_{\varepsilon}(r)h(r),\quad \bar{h}_{\varepsilon}(r):=\bar\phi_{\varepsilon}(r)h(r):=(1-\phi_{\varepsilon}(r))h(r), 
\]
where $\phi_{\varepsilon}\in C^{\infty}(\mbR)$ is a ``truncation" function such that $\mathbf{1}_{|r|\geqslant\varepsilon}\leqslant\phi_{\varepsilon}(r)\leqslant\mathbf{1}_{|r|\geqslant\varepsilon/2}$, $\phi_{\varepsilon}(w)$ is non-decreasing as $|w|$ increases, and ${\rm supp}(\phi_{\varepsilon})=(\varepsilon/2,\infty)\cup(-\infty,-\varepsilon/2)$. Next,
let $Z(\varepsilon):=\{Z_t(\varepsilon)\}_{t\geqslant0}$ and $Z'(\varepsilon):=\{Z'_t(\varepsilon)\}_{t\geqslant0}$ be independent pure-jump L\'{e}vy processes, with respective L\'{e}vy triplets
$(b_{Z}(\varepsilon),h_{\varepsilon}(r)\,dr,0)$ and $(0,\bar{h}_{\varepsilon}(r)\,dr,0)$, for the truncation function $\mathbf{1}_{|r|\leqslant1}$, where $b_{Z}(\varepsilon):=\int_{|r|\leqslant{}1}rh_{\varepsilon}(r)dr$. Note that $Z(\varepsilon)$ is a compound Poisson process, and we shall denote its intensity by $\lambda_{\varepsilon}:=\int\phi_{\varepsilon}(r)h(r)dr$, and the jump probability density function by $\breve{h}_{\varepsilon}(r):=\phi_{\varepsilon}(r)h(r)/\lambda_{\varepsilon}$.
Let $\{\tau_i\}_{i\geqslant1}$, $\{N_t\}_{t\geqslant0}$, and $\{J_i\}_{i\geqslant1}$
denote, respectively, the jump times, jump counting process of $Z(\varepsilon)$, and an independent identically distributed random sample from the probability density function $\breve{h}_{\varepsilon}$. Also,  $U$ and $J:=J^{\varepsilon}$ represent a generic random variable uniformly distributed in $(0,1)$ and a generic random variable with the probability density function $\breve{h}_{\varepsilon}(r)$, respectively.

The lemma below from \cite{Cheng} will be useful in the sequel.
\begin{lem}\label{gammainfinity}
Under the conditions (\ref{conditionnu}) and (\ref{conditiongamma}) in section \ref{setup}, the following statements hold:
\begin{enumerate}
\item Let $\tilde{\gamma}(z,r):=z+\gamma(z,r)$. Then, for each $z\in\mbR$, the mapping $r\to\tilde{\gamma}(z,r)$ (resp., $r\to\gamma(z,r)$) is invertible and its inverse $\tilde\gamma^{-1}(z,r)$ (resp., $\gamma^{-1}(z,r)$) belongs to $\bar{C}^{4}_{b}(\mbR\times\mbR_0)$.

\item Both $\tilde\gamma(z,J)$ and $\gamma(z,J)$ admit densities, denoted by $\widetilde\Gamma(r;z):=\widetilde\Gamma_{\varepsilon}(r;z)$ and $\Gamma(r;z):=\Gamma_{\varepsilon}(r;z)$, respectively, which belong to $C^{4}_{b}(\mbR\times\mbR_0)$. Furthermore, they have the representations:
\begin{equation*} \widetilde\Gamma_{\varepsilon}(r;z)={\breve{h}}_{\varepsilon}(\tilde{\gamma}^{-1}(z,r))\left|\frac{\partial \gamma}{\partial{}r}\left(z,\tilde\gamma^{-1}(z,r)\right)\right|^{-1},\quad\Gamma_{\varepsilon}(r;z)={\breve{h}}_{\varepsilon}({\gamma}^{-1}(z,r))\left|\frac{\partial \gamma}{\partial{}r}\left(z,\gamma^{-1}(z,r)\right)\right|^{-1}.
\end{equation*}
\item The mapping $z\to u:=z+\gamma(z,r)$ admits {an} {inverse}, denoted hereafter by $\bar{\gamma}(u,r)$, that belongs to $\bar{C}^{4}_{b}(\mbR\times\mbR_0)$.
	\end{enumerate}
\end{lem}

Now, we are ready to define the ``small-jump component" of $X^{(x)}$.
Let $M'(dv,dr):=M'_{\varepsilon}(dv,dr)$ denote the jump measure of the process $Z'(\varepsilon)$, and let $p'(dv,dr,du):=p'_{\varepsilon}(dv,dr,du)$ (resp. {$\bar{p}'(dv,dr,du):=p'_{\varepsilon}(dv,dr,du)-dv\bar{h}_{\varepsilon}(r)drdu$}) denote the marked point process (resp. compensated marked point process) on the atoms of $M'$ with independent uniformly distributed marks on $(0,1)$. For each $\varepsilon\in(0,1)$, we construct a process $X^{\theta}:=\left\{X_{s}^{\theta}(\varepsilon,x)\right\}_{s\geqslant0}$, defined as the solution of the SDE
\begin{equation}\label{Xtheta}
\begin{aligned}
X^{\theta}_{s}(\varepsilon,x)&=x+\int_{0}^{s} b_{\varepsilon}(X^{\theta}_{v^-}(\varepsilon,x))dv+ \int_{0}^{s}\sigma\left(X^{\theta}_{v^-}(\varepsilon,x)\right)d\widetilde{W}_{v}\\ &\quad+\sum_{i=1}^{N_{s}}\gamma\left(X^{\theta}_{\tau_{i}^{-}}(\varepsilon,x),J_{i}\right)\theta\left(X^{\theta}_{\tau_{i}^{-}}(\varepsilon,x),J_i,U_i\right)\\
&\quad+\int_{0}^{s}\int_E\gamma\left(X^{\theta}_{v^{-}}(\varepsilon,x),r\right)\theta\left(X^{\theta}_{v^{-}}(\varepsilon,x),r,u\right)\bar{p}'(dv,dr,du),
\end{aligned}
\end{equation}
where $\{\widetilde{W}_v\}_{v\geqslant0}$ is a Wiener process independent of $p'(dv,dr,du)$, and
\begin{equation}\label{bepsilon} b_{\varepsilon}(x)
:=b(x)-\int_{|r|\leqslant1}\int_0^1\gamma\left(x,r\right)\theta(x,r,u)duh_{\varepsilon}(r)dr
=b(x)-\int_{|r|\leqslant1}\gamma\left(x,r\right)\nu(x,r)\phi_{\varepsilon}(r)dr.
\end{equation}
By comparing their infinitesimal generators, it is not hard to see that  the process (\ref{Xtheta}) has the same distribution law as the process (\ref{SDEpq}) (see Section 2 of \cite{Cheng} for a more detailed explanation).
Next, we let $X^{\theta}(\varepsilon,\emptyset,x):=\left\{X^{\theta}_{s}(\varepsilon,\emptyset,x)\right\}_{s\geqslant{}0}$ be the solution of the SDE:
\begin{equation}\label{X0bigjump}
\begin{aligned}
X_{s}^{\theta}(\varepsilon,\emptyset,x)&=x+\int_{0}^{s}b_{\varepsilon}(X_{v^-}^{\theta}(\varepsilon,\emptyset,x))dv+\int_{0}^{s}\sigma\left(X_{v^-}^{\theta}(\varepsilon,\emptyset,x)\right)d\widetilde{W}_{v}\\ &\quad+\int_{0}^{s}\int_E\gamma\left(X_{v^{-}}^{\theta}(\varepsilon,\emptyset,x),r\right)\theta\left(X^{\theta}_{v^{-}}(\varepsilon,\emptyset,x),r,u\right)\bar{p}'(dv,dr,du).
\end{aligned}
\end{equation}

The law of the process $\left\{X_{s}^{\theta}(\varepsilon,\emptyset,x)\right\}_{0\leqslant{}s\leqslant{}t}$ above can be interpreted as  the law of $\{X^{\theta}_s(\varepsilon,x)\}_{0\leqslant{}s\leqslant{}t}$ conditioning on $\{Z_s(\varepsilon)\}_{0\leqslant{}s\leqslant{}t}$ not having any jumps. Note that, by conditions (\ref{Cinfinity}) and (\ref{conditiongamma}-\ref{gamma0=0}), the process (\ref{X0bigjump}) is a local martingale with bounded drift whose jumps are bounded by a constant.  With equation (9) in \cite{Lepeltier} followed by the proofs of Proposition 3.1 and Lemma 3.2 in \cite{Cheng}, we have
\begin{equation}\label{ttoN}
\sup_{0<\eta<\varepsilon,x\in\mbR}\mbP[|X_{t}^{\theta}(\eta,\emptyset,x)-x|\geqslant{}y]<C t^{N}
\end{equation}
for any $y>0$, where $N>0$ can be made arbitrarily large by taking $\varepsilon>0$ small enough.

We now proceed to define other related processes. For a collection of times $0<s_{1}<\dots<s_{n}$, let $\left\{X^{\theta}_{s}(\varepsilon,\{s_{1},\dots,s_{n}\},x)\right\}_{s\geqslant{}0}$ be the solution of the SDE:
\begin{equation*}
\begin{aligned}
X^{\theta}_{s}(\varepsilon,\{s_{1},\dots,s_{n}\},x)&=x+\int_{0}^{s}b_{\varepsilon}(X^{\theta}_{v^-}(\varepsilon,\{s_{1},\dots,s_{n}\},x))dv+\int_{0}^{s}\sigma\left(X^{\theta}_{v^-}(\varepsilon,\{s_{1},\dots,s_{n}\},x)\right)d\widetilde{W}_{v}\\ &\quad+\sum_{i:s_{i}\leqslant{}s}\gamma\left(X^{\theta}_{s_i^-}(\varepsilon,\{s_{1},\dots,s_{n}\},x),J_{i}\right)\theta\left(X^{\theta}_{s_i^-}(\varepsilon,x),J_i,U_i\right)\\
&\quad+\int_{0}^{s}\int_E\gamma\left(X^{\theta}_{v^{-}}(\varepsilon,\{s_{1},\dots,s_{n}\},x),r\right)\theta\left(X^{\theta}_{v^{-}}(\varepsilon,\{s_{1},\dots,s_{n}\},x),r,u\right)\bar{p}'(dv,dr,du).
\end{aligned}
\end{equation*}
The law of the process $\left\{X^{\theta}_{s}(\varepsilon,\{s_{1},\dots,s_{n}\},x)\right\}_{0\leqslant{}s\leqslant{}t}$ can be interpreted as  the law of $\{X^{\theta}_{s}(\varepsilon,x)\}_{0\leqslant{}s\leqslant{}t}$ conditioning on $\{Z_s(\varepsilon)\}_{0\leqslant{}s\leqslant{}t}$ having $n$ jumps at the times $0<s_1<s_2<\cdots<s_n<t$.

For future reference, let us remark that the infinitesimal generator of the small-jump component $\{X_t^{\theta}(\varepsilon,\emptyset,x)\}_{t\geqslant0}$, hereafter denoted by $L_{\varepsilon}$, can be written as
\begin{equation}\label{generator0bigjump}
\begin{aligned}
&L_{\varepsilon}f(y):=\mathcal{D}_{\varepsilon}f(y)+\mathcal{I}_{\varepsilon}f(y)\quad\text{with}\\
&\mathcal{D}_{\varepsilon}f(y)=b_{\varepsilon}(y)f'(y)+\frac{\sigma^2(y)}{2}f''(y),\\
&\mathcal{I}_{\varepsilon}f(y)=\int\left[f(y+\gamma(y,r))-f(y)-f'(y)\gamma(y,r)\right]\nu(y,r)\bar{\phi}_{\varepsilon}(r)dr,
\end{aligned}
\end{equation}
where $b_{\varepsilon}$ is defined in (\ref{bepsilon}) and $\bar{\phi}_{\varepsilon}(r)=1-\phi_{\varepsilon}(r)$.
Note that, for $f\in C^{2}_{b}(\mbR)$, $\mathcal{I}_{\varepsilon}f$ can be written as
\begin{align} 
\mathcal{I}_{\varepsilon}f(y)
&=\int\int_{0}^{1}\left[f''(y+\gamma(y,r\beta))(\partial_{2}\gamma(y,r\beta))^{2}+f'(y+\gamma(y,r\beta))\partial^{2}_{2}\gamma(y,r\beta)\right.\label{Irepresentation2}\\
&\;\quad\qquad\quad\left.-f'(y)\partial^{2}_{2}\gamma(y,r\beta)\right](1-\beta)d\beta\nu(y,r)\bar{\phi}_{\varepsilon}(r)r^2dr,\nonumber
\end{align}
which is finite and bounded due to the conditions (\ref{conditionnu}-\ref{hExistence}) and (\ref{conditiongamma}-\ref{gamma0=0}).

The following  first and second order Dynkin's formula for the process $\{X^{\theta}_t(\varepsilon,\emptyset,x)\}_{t\geqslant0}$ will be needed in the sequel:
\begin{align}
\label{Dynkin1}&\mbE\left[f(X^{\theta}_{t}(\varepsilon,\emptyset,x))\right]=f(x)+t\int_{0}^{1}\mbE\left[L_{\varepsilon}f(X^{\theta}_{\alpha{}t}(\varepsilon,\emptyset,x))\right]d\alpha,\quad\forall{}f\in{}C_b^2(\mbR),\\
\label{Dynkin2}&\mbE\left[f(X^{\theta}_{t}(\varepsilon,\emptyset,x))\right]=f(x)+tL_{\varepsilon}f(x)+t^{2}\int_{0}^{1}(1-\alpha)\mbE\left[(L_{\varepsilon})^2f(X^{\theta}_{\alpha{}t}(\varepsilon,\emptyset,x))\right]d\alpha,\quad\forall{}f\in{}C_b^4(\mbR).
\end{align}
Furthermore, for $f\in C_b^2(\mbR)$ (resp., $f\in{}C_b^4(\mbR)$), $\sup_{x}|L_{\varepsilon}f(x)|<\infty$ (resp., $\sup_{x}|L^{2}_{\varepsilon}f(x)|<\infty$) and, thus, the reminders in (\ref{Dynkin1}) (resp., (\ref{Dynkin2})) is $O(t)$ (resp., $O(t^{2})$) uniformly on $x$ and $t$. The proofs of (\ref{Dynkin1})-(\ref{Dynkin2}) follows from It\^o's formula and goes along the same lines as in the proof of Lemma 3.3 in \cite{Cheng}.

\begin{rem}
By the condition (\ref{conditiongamma}-\ref{partial2gamma>eta}) in Section \ref{setup}, the mapping $r\to\gamma(x,r)$ is a bijection in $\mbR$. Then, the generator (\ref{generator0bigjump}) can be rewritten as
\begin{align}\label{generator2}
L_{\varepsilon}f(y)=b_{\varepsilon}(y)f'(y)+\frac{\sigma^2(y)}{2}f''(y)+\int{}\left(f(y+z)-f(y)-zf'(y)\right)K_{\varepsilon}(y,dz),
\end{align}
where the kernel $K_{\varepsilon}(y,dz)$ can be written in the following two equivalent forms 
for any $A\in\mathcal{B}(\mbR)$,
\begin{align}\label{KAgamma}
 K_{\varepsilon}(y,A)=\int_{\mbR}\mathbf{1}_A(\gamma(y,r))\nu(y,r)\bar{\phi}_{\varepsilon}(r)dr	=\int_{A}\nu\left(y,\gamma^{-1}(y,z)\right)\bar{\phi}_{\varepsilon}\left(\gamma^{-1}(y,z)\right)\partial_z\gamma^{-1}(y,z)dz.
\end{align}
\end{rem}

\section{A weak solution process}\label{weaksol}
The main purpose of this section is to introduce an approach to overcome the difficulties posed by the discontinuous jump component $\gamma(x,r)\theta(x,r,u)$ in (\ref{Xtheta}) and (\ref{X0bigjump}). To this end, {we first show that the state-dependent local jump diffusions (\ref{Xtheta}) and (\ref{X0bigjump}) are equivalent to a state-free local jump-diffusion of the form (\ref{SDEIntroduction}) and prove some needed regularity on its coefficients.}

Note that the process $\{X_t^{\theta}(\varepsilon,\emptyset,x)\}_{t\geqslant0}$
is a semimartingale (see, e.g., \cite[III.2.18]{Jacod}) and, furthermore, comparing the generator (\ref{generator2}) to that in \cite[IX.4.6]{Jacod}, it is a homogeneous diffusion process with jumps as defined in \cite[IX.4.1]{Jacod}.
Then, we can determine the semimartingale characteristics $(B^*,C^*,\vartheta^*)$ of $\{X_t^{\theta}(\varepsilon,\emptyset,x)\}_{t\geqslant0}$,
relative to the identity truncation function, as
\begin{equation}\label{semicharacteristics} B_t^*=\int_0^tb_{\varepsilon}\left(X^{\theta}_s(\varepsilon,\emptyset,x)\right)ds,
\quad C_t^*=\int_0^t{}\sigma^2\left(X^{\theta}_s(\varepsilon,\emptyset,x)\right)ds,
\quad \vartheta^*(dt,dr)=K_{\varepsilon}\left(X^{\theta}_{t^{-}}(\varepsilon,\emptyset,x),dr\right)dt.
\end{equation}
By Definition III.2.24 and Theorem III.2.26 in \cite{Jacod}, a semimartingale with characteristics specified by (\ref{semicharacteristics}) is a weak solution of the SDE
\begin{equation}\label{Y0bigjump}
Y_t(\varepsilon,\emptyset,x)=\,x+\int_0^t
b_{\varepsilon}(Y_s(\varepsilon,\emptyset,x))ds+\int_0^t
\sigma(Y_s(\varepsilon,\emptyset,x))dW^{*}_s+\int_0^t\int{}\delta(Y_{s^-}(\varepsilon,\emptyset,x),r)\overline{\mu}^{*}(ds,dr),
\end{equation}
where, under the solution measure $\mbP^*$ in the canonical space $(\Omega^{*},\mathcal{F}^{*},\mbP^*)$,  $\{W^{*}_s\}_{s\geqslant0}$ is a one-dimensional Wiener process, $\mu^{*}(ds,dr)$ is an independent Poisson random measure on $\mbR_+\times\mbR_0$ with intensity measure $dsF(dr)$ and corresponding compensated measure $\overline{\mu}^{*}(ds,dr):=\mu^{*}(ds,dr)-dsF(dr)$. Here, $F$ is a positive $\sigma$-finite measure on $\mbR_{0}$ to be chosen below, while $\delta:\mbR\times\mbR_0\mapsto\mbR$ is a Borel function implicitly determined by $K$, $\delta$, and $F$ via the equation
\begin{align}\label{Kdelta}
   K_{\varepsilon}(x,A)=\int\mathbf{1}_{A}(\delta(x,r))F(dr),\quad\forall A\in\mathcal{B}(\mbR_0),
\end{align}
where we recalled from (\ref{KAgamma}) that $K_{\varepsilon}(y,A):=\int_{\mbR}\mathbf{1}_A(\gamma(y,r))\nu(y,r)\bar{\phi}_{\varepsilon}(r)dr$. In what follows, we take 
\[
	F(dr)=\bar{\phi}_{\varepsilon}(r)h(r)dr,
\]
for {$\varepsilon<\epsilon_{0}$, with $\epsilon_{0}$} given as in Condition (\ref{conditionh}).

In order to identify the function $\delta(x,r)$ corresponding to the above measure $F$, we introduce {the following functions $\psi:\mbR_0\to\mbR_0$ and $\bar{\psi}:\mbR\times\mbR_0\to\mbR_0$:}
\begin{equation}\label{DfnPhiBarPhi}
\psi(w)=\left\{
\begin{aligned}
&-\int_w^{\infty}\bar{\phi}_{\varepsilon}(r)h(r)dr,\quad{}w>0,\\
&\int_{-\infty}^w{}\bar{\phi}_{\varepsilon}(r)h(r)dr,\quad{}w<0,
\end{aligned}
\right.\qquad\quad
\bar{\psi}(x,w)=\left\{
\begin{aligned}
&-\int_{w}^{\infty}\bar{\phi}_{\varepsilon}(r)\nu(x,r)dr,\quad{}w>0,\\
&\int^{w}_{-\infty}\bar{\phi}_{\varepsilon}(r)\nu(x,r)dr,\quad{}w<0.
\end{aligned}
\right.
\end{equation}
Note that the
restrictions of the two mappings $w\mapsto\psi(w)$ and $w\mapsto\bar{\psi}(x,w)$ on $D_{\varepsilon}:=(-\varepsilon,0)\cup (0,\varepsilon)$ are strictly increasing and one-to-one, with range $\mbR_{0}$. Thus, they
admit ``local" inverses defined on $\mbR_0$ with range $D_{\varepsilon}$, hereafter, denoted by $\psi^{-1}(x,w)=\psi_{\varepsilon}^{-1}(x,w)$ and $\bar{\psi}^{-1}(x,w)=\bar{\psi}_{\varepsilon}^{-1}(x,w)$, respectively.

{From the form of the SDE (\ref{Y0bigjump}) and the fact that the support of the mean measure $ds\bar{\phi}_{\varepsilon}(r)h(r)dr$ of $\mu^*$ lies in $(-\varepsilon,\varepsilon)$, it is clear that the values of $\delta(x,r)$ for $|r|\geq{}\varepsilon$ are superfluous and we only need to define $\delta(x.r)$ for  $r\in D_{\varepsilon}$}. Let $\delta:\mbR\times{}D_{\varepsilon}\to\mbR$ be defined as
\begin{equation}\label{DefnDelta}
	\delta(x,w):=\gamma\left(x,\bar{\psi}^{-1}(x,\psi(w))\right).
\end{equation}
In order to show that the above function satisfies (\ref{Kdelta}), {observe} that, for each $x\in\mbR$ and $w\in{}D_{\varepsilon}$, the mapping $w\to\delta(x,w)$ is strictly monotone and, thus, its inverse, hereafter denoted by $\delta^{-1}(x,w)$, exists and satisfies
\begin{equation}\label{deltainverse}
\delta^{-1}(x,w)=\psi^{-1}\left(\bar{\psi}(x,\gamma^{-1}(x,w))\right).
\end{equation}
Then, $\psi(\delta^{-1}(x,w))=\bar{\psi}(x,\gamma^{-1}(x,w))$ and, from the definitions of $\psi$ and $\bar{\psi}$, we have the identity
\begin{equation*}
\left\{
\begin{aligned}
&-\int_{\delta^{-1}(x,w)}^{\infty}\bar{\phi}_{\varepsilon}(z)h(z)dz,\quad{}w>0,\\
&\int_{-\infty}^{\delta^{-1}(x,w)}\bar{\phi}_{\varepsilon}(z)h(z)dz,\quad{}w<0,
\end{aligned}
\right.=
\left\{
\begin{aligned}
&-\int_{\gamma^{-1}(x,w)}^{\infty}\bar{\phi}_{\varepsilon}(z)\nu(x,z)dz,\quad{}w>0,\\
&\int^{\gamma^{-1}(x,w)}_{-\infty}\bar{\phi}_{\varepsilon}(z)\nu(x,z)dz,\quad{}w<0.
\end{aligned}
\right.
\end{equation*}
Upon differentiation with respect to $w$, we get
\begin{align}\label{gammadelta}	\bar{\phi}_{\varepsilon}(\delta^{-1}(x,w))h(\delta^{-1}(x,w))\partial_2\delta^{-1}(x,w)=\bar{\phi}_{\varepsilon}(\gamma^{-1}(x,w))\nu(x,\gamma^{-1}(x,w))\partial_2\gamma^{-1}(x,w),
\end{align}
which implies (\ref{Kdelta}) with the chosen measure $F(dr)={\bar{\phi}_{\varepsilon}(r)h(r)dr}$.

We now proceed to show some needed regularity properties of the function $\delta:\mbR\times{}D_{\varepsilon}\to\mbR$, with which the almost sure existence of a stochastic flow of diffeomorphisms associated to the SDE (\ref{Y0bigjump}) below is guaranteed by virtue of results in \cite{Malliavin} (see Lemma \ref{Y0diffeomorphism} below). The proof of 
Proposition \ref{deltaproperty} 
is deferred to the Appendix 
\ref{pfdeltaproperty}.

\begin{pro}\label{deltaproperty}
  Under the conditions (\ref{conditionnu}), (\ref{conditiongamma}), and (\ref{conditionh}) in  Section \ref{setup}, the function $\delta$ defined in (\ref{DefnDelta}) can be continuously extended on $\mbR\times{}(-\varepsilon,\varepsilon)$ with $\delta(x,0):=0$, for any $\varepsilon>0$. Furthermore, 
  for $\varepsilon>0$ small enough,
  \begin{equation}\label{1+partialxdelta}
   {\rm (i)}\;\sup_{x\in\mbR,w\in D_{\varepsilon}}|\partial_{2}\delta(x,w)|<\infty, \quad {\rm (ii)}\;|\partial_{1}^{i}\delta(x,w)|\leq{}k|w|,\quad   {\rm (iii)}\;|1+\partial_1\delta(x,w)|>\eta,
  \end{equation}
 for any $0\leqslant{}i\leqslant{}2$, $w\in D_{\varepsilon}$, and some constant $\eta,k>0$, independent of $x$ and $w$.
\end{pro}

\begin{rem}\label{WhyS4}
The main purpose of the assumption
(\ref{conditionh}) in Section \ref{setup} is to simplify the verification of the regularity of $\delta$ as stated in Proposition \ref{deltaproperty}. However, it is important to note that all what follows, including the main Theorem \ref{tailexpansion}, hold true if $\nu$ and $h$ are such that the function $F(x,w)=\bar{\psi}^{-1}(x,\psi(w))$ satisfies the following conditions:
\begin{equation}\label{TchnCnd}
	\sup_{x\in\mbR,w\in D_{\varepsilon}}|\partial_{w}F(x,w)|<\infty, 
	\quad 
	\sup_{x\in\mbR,w\in D_{\varepsilon}}|\partial_{x}^{i}F(x,w)|<\infty, \quad i=0,1,2.
\end{equation}
In particular, if $\nu(x,r)=h(r)$, for an arbitrary  L\'evy density $h$ that is smooth enough outside any neighborhood of the origin, then (\ref{TchnCnd}) trivially holds true since $F(x,w)=w$. In that case, we recover the results in \cite{Cheng}.
\end{rem}

\begin{rem}\label{EILTP}
Using similar arguments to those at the beginning of this section, it is not hard to check that, under the condition (\ref{conditiongamma}-\ref{partial2gamma>eta}), the model (\ref{SDEpq}) is equivalent in law to the model (\ref{SDEIntroduction}) with $\gamma$ replaced by an appropriate function $\Upsilon$. Concretely, we need to take
\[
	\Upsilon(x,w)=\gamma\left(x,\widetilde{\psi}^{-1}\left(x,\widehat{\psi}(w)\right)\right),
\]
where
\begin{equation*}
\widehat\psi(w):=\left\{
\begin{aligned}
&-\int_w^{\infty}h(r)dr,\quad{}w>0,\\
&\int_{-\infty}^w{}h(r)dr,\quad{}w<0,
\end{aligned}
\right.\qquad\quad
\widetilde{\psi}(x,w):=\left\{
\begin{aligned}
&-\int_{w}^{\infty}\nu(x,r)dr,\quad{}w>0,\\
&\int^{w}_{-\infty}\nu(x,r)dr,\quad{}w<0.
\end{aligned}
\right.
\end{equation*}
However, the regularity of $\Upsilon$ is harder to study than that of the function $\delta$ in (\ref{DefnDelta}). Indeed, for instance, the first-order partial derivatives of the function
$\widetilde{F}(x,w):=\widetilde{\psi}^{-1}(x,\widehat{\psi}(w))$ are given by
\begin{align*}		
\partial_{2}\widetilde{F}(x,w)=\frac{{h}(\widetilde{F}(x,w))}{\nu(x,\widetilde{F}(x,w))},
\qquad\partial_{1}\widetilde{F}(x,w)=-\frac{(\partial_{1}\widetilde\psi)(x,\widetilde{F}(x,w))}{\nu(x,\widetilde{F}(x,w))},
\end{align*}
and the behavior of $h(w)$ and $\nu(x,w)$ as $w\to\pm\infty$ is now relevant as well.
\end{rem}

\section{Tail estimate for the one-big-jump process}\label{importantlemma}

In this section, we give an expansion for the tail probability of the process $X^{\theta}$ defined in (\ref{Xtheta}) conditioned on $N_{s}$ having only one jump. The following lemma (proof in Appendix \ref{pfcounterpartA.1}) is a counterpart of Lemma A.1 in \cite{Cheng}, even though the proof developed in the present article is new and much simpler. Below, we recall the notation $v(x)=\sigma^{2}(x)/2$ and set
\begin{align*} \nu_{\varepsilon}(x,r):=\nu(x,r)\phi_{\varepsilon}(r),\quad \bar\nu_{\varepsilon}(x,r):=\nu(x,r)(1-\phi_{\varepsilon}(r)).
\end{align*}
\begin{lem}\label{counterpartA.1}
{With the notation introduced in Section \ref{notation&definition}, let}
\begin{align}\label{DefHtz}
  H(t,z,q):=\mbE\left[\frac{\nu(z,J)}{h(J)}\left.\mbP\left[X_t^{\theta}(\varepsilon,\emptyset,v)\geqslant{}q\right]\right|_{v=z+\gamma(z,J)}\right].
\end{align}
Then, {under the conditions (\ref{conditionnu})-(\ref{conditionh}) in Section \ref{setup},}
\begin{align}\label{Exp1J}
    H(t,z,q)=H_0(z;q)+tH_1(z;q)+t^2\mathcal{R}^{(2)}_t(z,q),
\end{align}
 where,
\begin{equation*}
\begin{aligned}
   H_0(z;q):&=\frac{1}{\lambda_{\varepsilon}}\int^{\infty}_{\gamma^{-1}(z,q-z)}\nu_{\varepsilon}(z,r)dr,
   \quad{}H_1(z;q):=D(z;q)+I(z;q)\\
   D(z;q):&=\frac{1}{\lambda_{\varepsilon}}\Big(\big[b_{\varepsilon}(q)-v'(q)\big]\nu_{\varepsilon}\left(z,\tilde{\gamma}^{-1}(z,q)\right)\partial_2\tilde{\gamma}^{-1}(z,q)\\
&\qquad\quad-v(q)\big[\partial_2\nu_{\varepsilon}\left(z,\tilde{\gamma}^{-1}(z,q)\right)\left(\partial_2\tilde{\gamma}^{-1}(z,q)\right)^2
+\nu_{\varepsilon}\left(z,\tilde{\gamma}^{-1}(z,q)\right)\partial^2_2\tilde{\gamma}^{-1}(z,q)\big]\Big),\\
   I(z;q):&=\frac{1}{\lambda_{\varepsilon}}\int\bigg[\int^q_{\bar{\gamma}(q,r)}\nu_{\varepsilon}\left(z,\tilde{\gamma}^{-1}(z,\eta)\right)\partial_2\tilde{\gamma}^{-1}(z,\eta)\bar{\nu}_{\varepsilon}(\eta,r)d\eta\\
  &\qquad\qquad-\nu_{\varepsilon}\left(z,\tilde{\gamma}^{-1}(z,q)\right)\partial_2\tilde{\gamma}^{-1}(z,q)\gamma(q,r)\bar{\nu}_{\varepsilon}(q,r)\bigg]dr,
\end{aligned}
\end{equation*}
and, for $\varepsilon>0$ small enough,
\begin{equation}\label{boundremainder}
	\limsup_{t\to{}0} \sup_{z,q}\left|\mathcal{R}^{(2)}_{t}\left(z,q\right)\right|<\infty,
\qquad \sup_{z,q}|H_{1}(z;q)|<\infty.
\end{equation}
\end{lem}

The proof of Lemma \ref{counterpartA.1}, which can be found in Appendix \ref{pfcounterpartA.1}, builds on the following key lemma.
The proof can be found in Appendix \ref{pfY0diffeomorphism}.

\begin{lem}\label{Y0diffeomorphism}
  Under the conditions {(\ref{conditionnu})-(\ref{conditionh})} in Section \ref{setup}, the SDE (\ref{Y0bigjump}) admits a unique solution, which in turn implies the existence of a unique weak solution of the SDE (\ref{SDEpq}). Moreover, for any $t>0$, the mapping $x\mapsto{}Y_t(\varepsilon,\emptyset,x)$ is a 
  diffeomorphism on $\mbR$.
\end{lem}

\section{The second order expansion for the tail probability}\label{secExp}
In this section, we obtain a second order expansion, in a short time $t$, for the tail probability $\mbP\left[X_t^{(x)}\geqslant{}x+y\right]$ for any $y>0$. The idea is to exploit the equivalence in law of $X^{(x)}$ and the process $X^{\theta}$ defined in (\ref{Xtheta}), and to condition on the number of jumps of $Z(\varepsilon)$. Specifically, we have
\begin{equation}\label{conditional}
\mbP\left[X^{(x)}_{t}\geqslant{}x+y\right]=\mbP\left[X^{\theta}_{t}(\varepsilon,x)\geqslant{}x+y\right]=\sum_{n=0}^{\infty}\mbP\left[\left.X^{\theta}_{t}(\varepsilon,x)\geqslant{}x+y\right|N_{t}^{\varepsilon}=n\right]\frac{(\lambda_{\varepsilon}t)^{n}}{n!}e^{-\lambda_{\varepsilon}t}.
\end{equation}
In order to present the expansion, let us first recall the notation $\nu_{\varepsilon}(x,r):=\nu(x,r)\phi_{\varepsilon}(r)$ and $\bar\nu_{\varepsilon}(x,r):=\nu(x,r)\bar\phi_{\varepsilon}(r)$ as well as the functions $\tilde\gamma$, $\gamma^{-1}$, $\tilde\gamma^{-1}$, and $\bar\gamma$, introduced in Lemma \ref{gammainfinity}. The following theorem (proof in Appendix \ref{pftailexpansion}) states the main result of this article.
\begin{thrm}\label{tailexpansion}
Under the conditions {(\ref{conditionnu})-(\ref{conditionh}) in {Section \ref{setup}}}, the following asymptotic expansion holds, for any $y>0$,
\begin{equation}\label{2ndOrdExp2}
\mbP\left[X^{(x)}_t\geqslant{}x+y\right]=tP_1(x,y)+\frac{t^2}{2}P_2(x,y)+o(t^2),
\end{equation}
as $t\to{}0$, where the coefficients $P_1(x,y)$ and $P_2(x,y)$ admit the following representations, for $\varepsilon>0$ small-enough:
\begin{align*}
P_1(x,y):=\int_{\{r:\gamma(x,r)\geqslant{}y\}}\nu_{\varepsilon}(x,r)dr,\qquad
P_2(x,y):=\mathcal{D}_{\varepsilon}(x,y)+\mathcal{J}_{\varepsilon}(x,y),
\end{align*}
with
\begin{align*}
\mathcal{D}_{\varepsilon}(x,y)
&=b_{\varepsilon}(x)\left[-\nu_{\varepsilon}\left(x,\gamma^{-1}(x,y)\right)\left[\partial_1\gamma^{-1}(x,y)-\partial_2\gamma^{-1}(x,y)\right]
    +\int^{\infty}_{\gamma^{-1}(x,y)}\partial_1\nu_{\varepsilon}(x,r)dr\right]\\
    &\quad+\frac{\sigma^2(x)}{2}\bigg[-\partial_2\nu_{\varepsilon}\left(x,\gamma^{-1}(x,y)\right)\left(\partial_1\gamma^{-1}(x,y)-\partial_2\gamma^{-1}(x,y)\right)^{2}+\int^{\infty}_{\gamma^{-1}(x,y)}\partial_1^2\nu_{\varepsilon}(x,r)dr\\
    &\hspace{2cm}-\nu_{\varepsilon}\left(x,\gamma^{-1}(x,y)\right)\left[\partial_1^{2}\gamma^{-1}(x,y)-2\partial_{1}\partial_{2}\gamma^{-1}(x,y)+\partial_2^{2}\gamma^{-1}(x,y)\right]\bigg]\\
    &\quad+\Big(b_{\varepsilon}(x+y)-\sigma(x+y)\sigma'(x+y)\Big)\nu_{\varepsilon}\left(x,\gamma^{-1}(x,y)\right)\partial_2\gamma^{-1}(x,y)\\
    &\quad-\frac{\sigma^2(x+y)}{2}\Big(\partial_{2}\nu_{\varepsilon}\left(x,{\gamma}^{-1}(x,y)\right)\big[\partial_2{\gamma}^{-1}(x,y)\big]^{2}
    +\nu_{\varepsilon}\left(x,{\gamma}^{-1}(x,y)\right)\partial^{2}_2{\gamma}^{-1}(x,y)\Big),
\end{align*}
\begin{align*}
\mathcal{J}_{\varepsilon}(x,y)&=\int\bigg[\int^{\infty}_{\gamma^{-1}(x+\gamma(x,r),y-\gamma(x,r))}\nu_{\varepsilon}(x+\gamma(x,r),r_1)dr_1
  -\int_{\gamma^{-1}(x,y)}^{\infty}\nu_{\varepsilon}(x,r_1)dr_1\\
  \nonumber&\hspace{1cm}+\gamma(x,r)\bigg(\nu_{\varepsilon}\left(x,\gamma^{-1}(x,y)\right)\left[\partial_1\gamma^{-1}(x,y)-\partial_2\gamma^{-1}(x,y)\right]
  -\int^{\infty}_{\gamma^{-1}(x,y)}\partial_1\nu_{\varepsilon}(x,r_1)dr_1\bigg)\bigg]\bar{\nu}_{\varepsilon}(x,r)dr\\
  &\quad+\int\bigg[\int^{x+y}_{\bar{\gamma}(x+y,r)}\nu_{\varepsilon}\left(x,\tilde{\gamma}^{-1}(x,r_1)\right)\partial_2\tilde{\gamma}^{-1}(x,r_1)\bar\nu_{\varepsilon}(r_{1},r)dr_{1}\\
  &\quad\qquad\qquad-\nu_{\varepsilon}\left(x,\gamma^{-1}(x,y)\right)\partial_2\gamma^{-1}(x,y)\gamma(x+y,r)\bar\nu_{\varepsilon}(x+y,r)\bigg]dr\\
  &\quad+\int\nu_{\varepsilon}(x,r_1)\int^{\infty}_{\gamma^{-1}(x+\gamma(x,r_1),y-\gamma(x,r_1))}\nu_{\varepsilon}(x,r_2)dr_2dr_1\\
  &\quad-\int_{\gamma^{-1}(x,y)}^{\infty}\nu_{\varepsilon}(x,r_1)\int\nu_{\varepsilon}(\tilde\gamma(x,r_1),r_2)dr_2dr_1
  -\int^{\infty}_{\gamma^{-1}(x,y)}\nu_{\varepsilon}(x,r)dr\int\nu_{\varepsilon}(x,r)dr.
\end{align*}
\end{thrm}

\begin{rem}\label{SRemTP}\hfill
\begin{enumerate}
\item The expansion (\ref{2ndOrdExp2}) indeed does not depend on $h$, which is expected since the infinitesimal generator (\ref{InfGenX}) of $X$ only depends on $(b,\sigma,\gamma,\nu)$. Obviously, the coefficients $P_1$ and $P_{2}$ are independent of $\varepsilon$, even though, the given representations involve $\varepsilon$. In particular, note that, since $\gamma(x,0)=0$, for $\varepsilon>0$ small-enough, $P_1(x,y)=\int_{\{r:\gamma(x,r)\geqslant{}y\}}\nu(x,r)dr.$
\item The drift and diffusion are absent in the leading term, which can be interpreted by saying that a possible ``big" move of the process in a short time happens mostly as a result of a ``large" jump. This phenomenon also appears in the expansion of \cite{Cheng}.
In the case that $\nu(x,r)$ is interpreted as the probability density of a mark $J$ of the underlying marked point process driving $X$, then $P_{1}(x,y)=\mbP\left[\gamma(x,J)\geqslant{}y\right]$.

\item It is tedious but not hard to verify that the above expansion reduces to the expansion of \cite{Cheng} when $\nu(x,r)=h(r)$, in which case the jump intensity does not depend on the state of the process $\{X_t\}_{t\geqslant0}$.
\item When $\gamma(x,r)=r$ and $\varepsilon\in(0,1)$ is small enough, the expansion reduces to
\begin{align*}
&t\int_{y}^{\infty}\nu(x,r)dr\\
&+\frac{t^2}{2}\bigg\{b_{\varepsilon}(x)\bigg[\nu\left(x,y\right)
    +\int^{\infty}_{y}\partial_1\nu(x,r)dr\bigg]+\frac{\sigma^2(x)}{2}\bigg[-\partial_2\nu\left(x,y\right)+\int^{\infty}_{y}\partial_1^2\nu(x,r)dr\bigg]\\
    &\qquad\quad+\Big(b_{\varepsilon}(x+y)-\sigma(x+y)\sigma'(x+y)\Big)\nu\left(x,y\right)-\frac{\sigma^2(x+y)}{2}\partial_{2}\nu\left(x,y\right)
\\
&\qquad\quad+\int\bigg[\int^{\infty}_{y-r}\nu(x+r,r_1)dr_1
  -\int_{y}^{\infty}\nu(x,r_1)dr_1-r\bigg(\nu\left(x,y\right)
  +\int^{\infty}_{y}\partial_1\nu(x,r_1)dr_1\bigg)\bigg]\bar{\nu}_{\varepsilon}(x,r)dr\\
  &\qquad\quad+\int\bigg[\int^{y}_{y-r}\nu_{\varepsilon}\left(x,r_{1}\right)\bar\nu_{\varepsilon}(x+r_{1},r)dr_{1}-\nu_{\varepsilon}\left(x,y\right)r\bar\nu_{\varepsilon}(x+y,r)\bigg]dr\\
  &\qquad\quad+\int\nu_{\varepsilon}(x,r_1)\int^{\infty}_{y-r_1}\nu_{\varepsilon}(x,r_2)dr_2dr_1\\
  &\qquad\quad-\int_{y}^{\infty}\nu(x,r_1)\int{\nu}_{\varepsilon}(x+r_1,r_2)dr_2dr_1
  -\int^{\infty}_{y}\nu(x,r)dr\int{\nu}_{\varepsilon}(x,r)dr\bigg\}+o(t^2).
\end{align*}
In particular, supposed that $b$ and $\sigma$ are constants. Then, recalling that $b_{\varepsilon}=b-\int_{|r|\leqslant1}r\nu(x,r)\phi_{\varepsilon}(r)dr$, the effect of a positive ``drift" $b$ is to increase the probability of a ``large" move of more than $y$ by
\[
	\frac{t^{2}}{2}b\left(2\nu(x,y)+\int_{y}^{\infty}\frac{\partial}{\partial x}\nu(x,r)dr\right)\left(1+o(1)\right).
\]
Note that the second term inside the parentheses above is missing in the absence of a state-dependent jump intensity. Similarly,
a nonzero constant volatility $\sigma$ will change the probability of a ``large" move of more than $y$ by
\begin{align*} \frac{t^{2}}{2}\sigma^{2}\left(-\frac{\partial}{\partial y}\nu(x,y)+\frac{1}{2}\int_{y}^{\infty}\frac{\partial^{2}}{\partial x^{2}}\nu(x,r)dr\right)\left(1+o(1)\right),
\end{align*}
in short-time. Again, the second term inside the parenthesis is the effect of a state-dependent intensity. For a general $b$, it is also intuitive that the net effect depends on the average drift, $(b(x)+b(x+y))/2$, at the initial and final points $x$ and $x+y$. The net effect of a general function $\sigma$ on the probability of a large positive move of more than $y$ depends on the functional
\[
\Lambda(x,y):=-\frac{\partial}{\partial{}y}\left(\nu(x,y)\,\frac{\sigma^{2}(x)+\sigma^{2}(x+y)}{2}\right)
+\frac{\sigma^{2}(x)}{2}\int_{y}^{\infty}\frac{\partial^{2}}{\partial{}x^{2}}\nu(x,r)dr,
\]
as $\frac{t^{2}}{2}\Lambda(x,y)\left(1+o(1)\right)$.
\end{enumerate}
\end{rem}

\section{The small-time second-order expansion for OTM call prices}\label{option}
In this section, we derive a second-order expansion, in short-time, for the price of an out-of-the-money (OTM) European call option, with maturity $t$ and strike $K$, written on a nondividend paying stock, whose risk-neutral price process is modeled by
\begin{align*}
S_t=S_{0}e^{X_t^{(0)}},\qquad t\geqslant0,
\end{align*}
where the process $\{X_t^{(0)}\}_{t\geqslant0}$ is given by (\ref{SDEpq}) with the initial condition $x=0$. For simplicity, in the rest of this section, we omit the superscript in $X^{(0)}$.

As {explained in the introduction}, we shall consider the share measure associated with the stock (namely, the martingale measure obtained by taking the stock as num\'eraire) {to evaluate the premium of the OTM option}. Concretely, let $k:=\log(K/S_0)$ be the so-called log-moneymess of this call option and as customary suppose that the risk-free rate $r$ is $0$. Then, the price of this option can be written as
\begin{align}\label{EsyRelOpPr}
\mathbb{E}[(S_t-K)_+]=
\mathbb{E}[(S_t-K)\mathbf{1}_{S_t\geqslant K}]
=S_0\mathbb{P}^{\#}\left[X_t\geqslant k\right]-S_0e^k\mathbb{P}\left[X_t\geqslant k\right],
\end{align}
where $\mathbb{P}^{\#}$ is a probability measure, locally equivalent to $\mathbb{P}$, defined as $\mathbb{P}^{\#}[B]:=\mathbb{E}[e^{X_t}\mathbf{1}_B]$  for any $\mathcal{F}_t$-measurable set $B$. Hereafter, $E^{\#}$ denotes the corresponding expectation. For $\mathbb{P}^{\#}$ to be well-defined, $\{S_{t}/S_{0}\}_{t\geqslant0}$ must be a $\mathbb{F}$-martingale, for which we impose the following {drift} restriction:
\begin{equation}\label{martingalecondition}
b(x)+\frac{\sigma^2(x)}{2}
+\int\left(e^{\gamma(x,r)}-1-\mathbf{1}_{|r|\leqslant1}\gamma(x,r)\right)\nu(x,r)dr=0.
\end{equation}
The integral in (\ref{martingalecondition}) will be well defined under the conditions (\ref{conditiongamma}-\ref{gamma0=0}) and (\ref{conditionnu}) in Section \ref{setup} as well as the condition:
\begin{align*}
	\int_{|r|\geqslant1}e^{\gamma(x,r)}\nu(x,r)dr<\infty.
\end{align*}
That $\{X_{t}\}_{t\geqslant0}$ is a $\mathbb{F}$-martingale under (\ref{martingalecondition}) is a consequence of It\^o's formula.

The second-order expansion for the probability appearing on the second term of (\ref{EsyRelOpPr}) was treated in Section \ref{setup}.  Next, we seek for a second order expansion in $t$ for the tail probability $\mathbb{P}^{\#}[X_t\geqslant k]$. To this end, we impose the following condition:
\renewcommand{\labelenumi}{(\theenumi)}
\renewcommand{\theenumi}{S\arabic{enumi}}
\renewcommand{\labelenumii}{(\theenumii)}
\renewcommand{\theenumii}{\roman{enumii}}
\begin{enumerate}
\setcounter{enumi}{4}
\item\label{expgamma} The function $g$ {introduced in Assumption (\ref{conditionh})} is such that $\int_{r\geqslant1}e^{cr}g(r)dr<\infty$,
where $c$ is defined by $c:=\sup\limits_{x,r}|\partial_2\gamma(x,r)|<\infty$.
    \renewcommand{\theenumi}{S\arabic{enumi}}
\end{enumerate}
\renewcommand{\labelenumi}{\theenumi.}
\renewcommand{\theenumi}{\arabic{enumi}}

Our first task is to determine the infinitesimal generator of the  process $\left\{ X_t\right\}_{t\geqslant{}0}$ under $\mathbb{P}^{\#}$. To that end, we compute the expectation
\begin{align*}
	\mbE^{\#}\left[q\left( X_t\right)\right]=\mbE\left[e^{{ X_t}}q\left(X_t\right)\right]
\end{align*}
for an arbitrary function $q\in C_b^2(\mbR)$. Let $f(x)=e^xq(x)$. Then, applying It\^{o}'s formula (see, e.g., \cite{Applebaum} Theorem 4.4.7),
$f\left({X_t}\right)=M_t+A_t$, where
\begin{align*}
M_t&=q(0)+\int_0^te^{X_s}\left[q\left(X_{s}\right)
+q'\left(X_s\right)\right]\sigma\left(X_s\right)dW_s\\
&\quad+\int_{[0,t)\times{}E}e^{X_{s^-}}\left[e^{\gamma\left(X_{s^-},r\right)}q\left(X_{s^-}+\gamma\left(X_{s^-},r\right)\right)-q\left(X_{s^-}\right)\right]\mathbf{1}_{\theta\left(X_{s^{-}},r,u\right)=1}\bar{p}(ds,dr,du),\\
A_t&=\int_0^t{}e^{X_{s}}\left[q\left(X_{s}\right)+q'\left(X_{s}\right)\right]b\left(X_{s}\right)ds
+\frac{1}{2}\int_0^t{}e^{X_{s}}\left[q\left(X_s\right)+2q'\left(X_{s}\right)+q''\left(X_{s}\right)\right]\sigma^2\left(X_{s}\right)ds\\
&\quad+\int_0^te^{X_{s}}\int\bigg(e^{\gamma\left(X_{s},r\right)}q\left(X_{s}+\gamma\left(X_{s},r\right)\right)-q\left(X_{s}\right)-\mathbf{1}_{|r|\leqslant1}\gamma\left(X_{s},r\right)\left[q\left(X_{s}\right)+q'\left(X_{s}\right)\right]\bigg)\nu\left(X_{s},r\right)drds.
\end{align*}
Thus,
\begin{align*}
\mbE^{\#}\left[q\left(X_t\right)\right]
&=q(0)+\int_0^t\mbE^{\#}\left[L^{\#}q\left(X_s\right)\right]ds,
\end{align*}
where
\begin{align*}
L^{\#}q(x)=\,&b(x)\left[q(x)+q'(x)\right]
+\frac{\sigma^2(x)}{2}\left[q(x)+2q'(x)+q''(x)\right]\\
&+\int\left[e^{\gamma(x,r)}q(x+\gamma(x,r))-q(x)-\mathbf{1}_{|r|\leqslant1}\gamma(x,r)[q(x)+q'(x)]\right]\nu(x,r)dr.
\end{align*}
Comparing the above formula to the Dynkin's formula (\ref{Dynkin1}), we can then identify $L^{\#}$ as the infinitesimal generator of $\{X_t\}_{t\geqslant0}$ under $\mbP^{\#}$. Using the martingale condition (\ref{martingalecondition}), we can further write $L^{\#}$ as
\begin{align}\label{generatorr}
L^{\#}q(x)
&=b^{\#}(x)q'(x)
+\frac{\sigma^2(x)}{2}q''(x)
+\int\left[q(x+\gamma(x,r))-q(x)-\mathbf{1}_{|r|\leqslant1}\gamma(x,r)q'(x)  \right]\nu^{\#}(x,r)dr,
\end{align}
where
\begin{align}
\label{bbound}
\nu^{\#}(x,r):=e^{\gamma(x,r)}\nu(x,r),\quad b^{\#}(x)
:=b(x)+\sigma^2(x)
+\int\left(e^{\gamma(x,r)}-1\right)\mathbf{1}_{|r|\leqslant1}\gamma(x,r)\nu(x,r)dr.
\end{align}
Note that under the {Conditions (\ref{conditiongamma}-\ref{gamma0=0}) and (\ref{conditionnu})} in Section \ref{setup}, the integral appearing in (\ref{bbound}) is well-defined and, furthermore, it is not hard to see that
$b^{\#}$ belongs to $C_b^4$ and can be written as
\begin{align}\label{bbound2}
b^{\#}(x)&=\frac{\sigma^2(x)}{2}
-\int\left(e^{\gamma(x,r)}-1-e^{\gamma(x,r)}\mathbf{1}_{|r|\leqslant1}\gamma(x,r)\right)\nu(x,r)dr,
\end{align}
in light of (\ref{martingalecondition}).

In order to conclude the second-order expansion for the call option price, it remains to show that the measure $\nu^{\#}$ satisfies the Conditions (\ref{conditionnu}) and (\ref{conditionh}). This is obtained in the following result whose proof is deferred to Appendix \ref{ProofPpoundexpansion}:
\begin{coro}\label{Ppoundexpansion}
Under the Conditions {(\ref{conditionnu})-(\ref{expgamma})}, for any $y>0$, we have
\begin{equation*}
\mbP^{\#}\left[X_t\geqslant{}y\right]=tP^{\#}_1(0,y)+\frac{t^2}{2}P_2^{\#}(0,y)+o(t^2),
\end{equation*}
as $t\to0$, where, $P_{1}^{\#}$ and $P_{2}^{\#}$ are given as Theorem \ref{tailexpansion}, but replacing $\nu$ and $b$ by $\nu^{\#}$ and $b^{\#}$, respectively, as defined in (\ref{bbound}).
\end{coro}
Finally, using the pricing formula for the European call option introduced at the beginning of this section, the price of an OTM European call option has the following expansion in a short maturity $t$, for $k>0$:
\begin{align}\nonumber
\mathbb{E}\left(S_t-S_{0}e^{k}\right)_+
&=S_0\mbP^{\#}\left[X_t\geqslant{}k\right]-S_0e^k\mathbb{P}\left[X_t\geqslant{}k\right]\\
\label{SecndOptPr}&=tS_0\big[P^{\#}_1(0,k)-e^kP_1(0,k)\big]
+\frac{t^2}{2}S_0\big[P^{\#}_2(0,k)-e^kP_2(0,k)\big]+o\big(t^{2}\big).
\end{align}
The formula (\ref{SecndOptPr}) extends the first order expansion for the option price given in \cite{Cheng} for a non state-dependent jump intensity process (i.e., when $\nu(x,r)=h(r)$).

\begin{rem}\label{remarkonS}\hfill
\begin{enumerate}
\item
It is not surprising that the leading term is only determined by the jump component according to the formula
\begin{align}\label{LOTCOP}
tS_{0}\left(\int_{\{r:\gamma(0,r)\geqslant{}k\}}e^{\gamma(0,r)}\nu(0,r)dr-e^{k}\int_{\{r:\gamma(0,r)\geqslant{}k\}}\nu(0,r)dr\right)
=tS_0\int\left(e^{\gamma(0,r)}-e^k\right)_+\nu(0,r)dr.
\end{align}
In particular, if $C(t,K):=\mathbb{E}\left(S_{t}-K\right)_+$ represents the premium of a call with expiration $t$ and strike $K$ and $\gamma(x,r)=r$, (\ref{SecndOptPr})-(\ref{LOTCOP}) suggest that, for $\kappa:=\log(K/S_0)>0$,
\begin{align*}
	\frac{\partial^{2}C(t,K)}{\partial{}K^{2}}\approx{}te^{-\kappa}\nu(0,\kappa) \quad\Longleftrightarrow
\quad\nu(0,\kappa)\approx\frac{1}{t}e^{\kappa}\frac{\partial^{2}C(t,K)}{\partial{}K^{2}};
\end{align*}
thus, the curvature of the call premium $K\to C(t,K)$ is strongly determined by jump intensity $\nu(0,\kappa)$.
\item Plugging the expansion given in Remark \ref{SRemTP}-4 into (\ref{SecndOptPr}) and using (\ref{martingalecondition}) and (\ref{bbound}), we note that, when $\gamma(x,r)={}r$, the effect of a nonzero constant volatility $\sigma(x)\equiv \sigma$ in the price of an OTM call option is of order
\begin{align*}	
{S_0\frac{t^2\sigma^2}{2}\bigg[e^k\nu(0,k)+\int_{k}^{\infty}\frac{e^r+e^k}{2}\partial_1\nu(0,r)dr+\int_{k}^{\infty}\frac{e^r-e^k}{2}\partial_1^2\nu(0,r)dr\bigg](1+o(1))}, \quad t\to0,
\end{align*}
which extends a result in \cite{Martin}  for exponential L\'evy models. We again observe an additional contribution due to the state dependent feature of the model.
\end{enumerate}
\end{rem}

\section{A numerical example of tail probability estimates}\label{example}
In this section, we show a comparison of numerical estimates for the tail probability $\mbP[X_t\geqslant{}y]$ by the first and the second order expansion in a small $t$, and a Monte Carlo estimate based on the simulation of $\{X_t\}_{t\geqslant0}$ with a jump augmented Euler-Maruyama scheme (cf. \cite{Mordecki}), combined with a diffusion approximation of the small-jump component of $X$. For the numerical results, we use the following parameters
\begin{equation}\label{1stparameters}
\begin{aligned}
&b(x)=\sin{}x,\quad\sigma(x)=\frac{1}{2}+\frac{1}{4}\sin{}x,\quad\gamma(x,r)=r,\\
&h(r)=|r|^{-1-\alpha}\,(\alpha=1.01),\quad{}\nu(x,r)=\left(\frac{1}{2\pi}tan^{-1}(x)+\frac{3}{4}\right)|r|^{-1-\alpha}=:c(x)h(r),
\end{aligned}
\end{equation}
which satisfy the conditions (\ref{conditionnu})-(\ref{conditionh}) in Section \ref{setup}. For simplicity, we set  $\phi_\varepsilon(r)=\mathbf{1}_{\{r:|r|>\varepsilon\}}$, when computing the coefficients $P_{1}(x,y)$ and $P_{2}(x,y)$ in the expansion (\ref{2ndOrdExp2}). This is valid since $P_{1}(x,y)$ and $P_{2}(x,y)$ don't depend on $\phi_{\varepsilon}$ and, thus, one {can} consider a sequence of smooth truncation functions, $\phi_{\varepsilon,n}$, that converges to $\phi_\varepsilon(r)=\mathbf{1}_{\{r:|r|>\varepsilon\}}$, as $n\to\infty$.

\subsection{Second order approximation for the tail probability}
With the parameters (\ref{1stparameters}) above and fixing $x=0$, we can compute explicitly the expansion stated in Remark \ref{SRemTP}-4, in which
\begin{align*}
  &\tilde\gamma(0,r)=\gamma(0,r)=r,\quad\tilde\gamma^{-1}(0,y)=\gamma^{-1}(0,y)=y,\quad\nu_{\varepsilon}(0,\gamma^{-1}(0,y))=\frac{3}{4}h_{\varepsilon}(y),\\
  &\partial_1\gamma(0,y)=\partial_1\gamma^{-1}(0,y)=\partial_1^i\partial_2^j\gamma(0,y)=\partial_1^i\partial_2^j\gamma^{-1}(0,y)=0 \quad(i+j=2),\\
  &\partial_2\gamma(0,y)=\partial_2\gamma^{-1}(0,y)=\partial_2\tilde\gamma^{-1}(0,y)=1,\quad\nu_{\varepsilon}(\tilde\gamma(0,r_1),r_2)=c(r_1)h_{\varepsilon}(r_2),\quad\partial_2\nu_{\varepsilon}(0,\gamma^{-1}(0,y))=\frac{3}{4}h'_{\varepsilon}(y).
\end{align*}
Then,
\begin{align*}
  P_1(0,y)&=\frac{3}{4}\int_{y}^{\infty}h_{\varepsilon}(r)dr,\\
  \mathcal{D}_{\varepsilon}(0,y)&=-\frac{3}{32}h'_{\varepsilon}(y)
  +\frac{3}{4}h_{\varepsilon}(y)[b_{\varepsilon}(y)-\sigma(y)\sigma'(y)]
  -\frac{3}{8}\sigma^2(y)h'_{\varepsilon}(y),\\
  \mathcal{J}_{\varepsilon}(0,y)&=\frac{3}{4}\int\bigg[
  c(r)\int^{\infty}_{y-r}h_{\varepsilon}(r_1)dr_1
  -\frac{3}{4}\int^{\infty}_{y}h_{\varepsilon}(r_1)dr_1
  -r\left(\frac{3}{4}h_{\varepsilon}(y)+\frac{1}{2\pi}\int^{\infty}_yh_{\varepsilon}(r_1)dr_1\right)\bigg]\bar{h}_{\varepsilon}(r)dr\\
  &\quad+\frac{3}{4}\int\bigg[\int_{y-r}^yc(r_1)h_{\varepsilon}(r_1)dr_1
  -c(y)h_{\varepsilon}(y)r\bigg]\bar{h}_{\varepsilon}(r)dr\\
  &\quad+\frac{9}{16}\int{}h_{\varepsilon}(r)\int^{\infty}_{y-r}h_{\varepsilon}(r_1)dr_1dr
  -\frac{3}{4}\lambda_{\varepsilon}\int_y^{\infty}h_{\varepsilon}(r)c(r)dr
  -\frac{9}{16}\lambda_{\varepsilon}\int_y^{\infty}h_{\varepsilon}(r)dr.
\end{align*}

\subsection{A Monte Carlo estimate of the tail probability}
In the sequel, we present a numerical method to simulate the process $X:=X^{(x)}$ defined in SDE (\ref{SDEpq}). This is based on a diffusion type approximation of the ``small-jump" component of the process together with a jump augmented Euler-Maruyama scheme (cf. \cite{Mordecki}) for the ``big-jump" component. To introduce the main idea, let us start by writing the infinitesimal generator $L$ of $X$, defined in (\ref{InfGenX}), as
\begin{align}
\nonumber{}L f(x)&=\hat{b}(x)f'(x)
+\frac{\sigma^2(x)}{2}f''(x)\\
\nonumber
&\quad +\int\Big[f(x+\gamma(x,r))-f(x)-\gamma(x,r)f'(x){\bf 1}_{\{|\gamma(x,r)|\leqslant{}1\}}\Big]\nu(x,r){\bf 1}_{\{|\gamma(x,r)|>{}\varepsilon\}}dr\\
\label{smallJGenerator}&\quad+\int\Big[f(x+\gamma(x,r))-f(x)-f'(x)\gamma(x,r)\Big]\nu(x,r){\bf 1}_{\{|\gamma(x,r)|\leqslant{}\varepsilon\}}dr,
\end{align}
for any $f\in{}C_b^2(\mbR)$ and $\varepsilon\in(0,1)$, where
	$\hat{b}(x):=b(x)-\int\gamma(x,r)[{\bf 1}_{\{|r|\leqslant{}1\}}-{\bf 1}_{\{|\gamma(x,r)|\leqslant{}1\}}]\nu(x,r)dr$.
Note that, as $\varepsilon\to0$, the term in (\ref{smallJGenerator}) can be approximated by
\begin{align*}
\frac{1}{2}f''(x)\int\gamma^2(x,r)\nu(x,r){\bf 1}_{\{|\gamma(x,r)|\leqslant{}\varepsilon\}}dr
  =:\frac{f''(x)}{2}\hat{\sigma}_{\varepsilon}^{2}(x).
\end{align*}
Therefore, for $\varepsilon$ small, the generator $L$ is ``close" to an operator $\tilde{L}_{\varepsilon}$ on $C_b^2(\mbR)$, defined as
\begin{equation}\label{InfGenTildeL}
\tilde{L}_{\varepsilon}f(x):= \hat{b}(x)f'(x)+\frac{\tilde\sigma^2_{\varepsilon}(x)}{2}f''(x)
+\int\Big[f(x+z)-f(x)-zf'(x){\bf 1}_{\{|z|\leqslant{}1\}}\Big]\widetilde{K}_{\varepsilon}(x,dz),
\end{equation}
where
\[ \tilde\sigma_{\varepsilon}^2(x):=\sigma^2(x)+\hat\sigma^{2}_{\varepsilon}(x),\quad
	\widetilde{K}_{\varepsilon}(x,A):=\int {\bf 1}_{A}(\gamma(x,r))\nu(x,r){\bf 1}_{\{|\gamma(x,r)|>{}\varepsilon\}}dr.
\]
The operator $\tilde{L}_{\varepsilon}$ is the infinitesimal generator of a Markov process $\{\tilde{X}_t(\varepsilon,x)\}_{t\geqslant0}$ of finite jump activity. Indeed, as stated in Remark \ref{remarkonS}, condition (\ref{conditiongamma}-\ref{gamma0=0}) implies that there exists an $\varepsilon_{0}>0$, depending only on $\varepsilon$ and not on $x$, such that $\{r:|\gamma(x,r)|>\varepsilon\}\subset\{r:|r|>\varepsilon_{0}\}$, for all $x\in\mathbb{R}$. Therefore, $\{\widetilde{X}_t^{\varepsilon}\}_{t\geqslant0}$ can be defined as the solution of the SDE:
\begin{align}\label{Xhat}
	\widetilde{X}_{t}^{\varepsilon}&=x+\int_{0}^{t} \tilde{b}_{\varepsilon}(\widetilde{X}_{v^-}^{\varepsilon})dv+ \int_{0}^{t}\tilde\sigma_{\varepsilon}\left(\widetilde{X}_{v^-}^{\varepsilon}\right)d\widetilde{W}_{v}+\sum_{i=1}^{\widetilde{N}^{\varepsilon}_{t}}\tilde{\gamma}_{\varepsilon}\left(\widetilde{X}_{\tau_{i}^{-}}^{\varepsilon},\widetilde{J}_{i}\right)\theta\left(\widetilde{X}_{\tau_{i}^{-}}^{\varepsilon},\widetilde{J}_i,\widetilde{U}_i\right),
\end{align}
where
	$\tilde{b}_{\varepsilon}(x)
		:=b(x)-\int\gamma(x,r)[{\bf 1}_{\{|r|\leqslant{}1\}}-{\bf 1}_{\{|\gamma(x,r)|\leqslant{}\varepsilon\}}]\nu(x,r)dr$,
$\tilde{\gamma}_{\varepsilon}(x,r):=\gamma(x,r){\bf 1}_{\{|\gamma(x,r)|>\varepsilon\}}$, $\{\widetilde{W_v}\}_{v\geqslant0}$ is a Wiener process, $\{\widetilde{N}_t^{\varepsilon}\}_{t\geqslant0}$ is an independent Poisson process with jump intensity $\lambda_{\varepsilon}=\int{}h(r){\bf 1}_{\{|r|>\varepsilon\}}dr$ and jump arrival times $\{\tau_{i}\}_{i\geqslant{}1}$, $\{\widetilde{J}_i\}_{i\geqslant1}$ are independent identically distributed with the probability density function $\breve{h}_{\varepsilon}(r)=h(r){\bf 1}_{\{|r|>\varepsilon\}}/\lambda_{\varepsilon}$, and $\{\widetilde{U}_i\}_{i\geqslant1}$ is a random sample from a standard uniform distribution. The following result, whose proof is deferred to Appendix \ref{pfsmallJ}, justifies that
$\{X_t\}_{t\geqslant0}$ can be asymptotically approximated in law by $\{\widetilde{X}_t^{\varepsilon}\}_{t\geqslant0}$, as $\varepsilon\to0$.
\begin{lem}\label{smallJ}
  Under the conditions (\ref{Cinfinity})-(\ref{conditionh}) in Section \ref{setup}, $\{\widetilde{X}_t^{\varepsilon}\}_{t\geqslant0}\xrightarrow{\mathcal{D}}\{X_t\}_{t\geqslant0}$.
\end{lem}

\subsection{Implementation and numerical results}
By Lemma \ref{smallJ}, we can approximate $\{X_t\}_{t\geqslant0}$ by simulating the process $\{\widetilde{X}_t\}_{t\geqslant0}$ defined in (\ref{Xhat}).
Since $\{\tilde{N}_t^{\varepsilon}\}_{t\geqslant0}$ is a Poisson process with intensity $\lambda_{\varepsilon}$, the inter-arrival jump times of (\ref{Xhat}) have independent identical exponential distribution with mean $1/\lambda_{\varepsilon}$. Then, we generate the jump times $\{\tau_{1},\tau_{2},\cdots,\tau_{m}\}$ by setting
\begin{align*}
  \tau_{0}=0,\quad
  \tau_{k}=\tau_{k-1}+\frac{1}{\lambda_{\varepsilon}}e_k,
\end{align*}
where $m:=\max\{k:\tau_k\leqslant{}t\}$ and $\{e_k\}_{k>0}$ is a random sample from the exponential distribution with parameter $1$. Also, we use the inverse transformation sampling method to obtain a sample of the jumps $\{\widetilde{J}_{i}\}_{i=1}^{m}$. Finally, we construct $\{\widetilde{X}_{t_k}\}_{k=1}^{n+m}$ over the jump-augmented time steps
\begin{align*}
  \{t_k\}_{k=1}^{n+m}:=\{s_j\}_{j=0}^n\cup\{\tau_i\}_{i=1}^m,\quad\left(s_j=j\frac{t}{n},\,n\in\mathbb{N}^+\right).
\end{align*}
The algorithm is similar to the one in \cite{Mordecki} but with an additional thinning condition at each jump time $\tau_i$.

\begin{figure}[H]
\centering
\includegraphics[width=5in,height=3in]{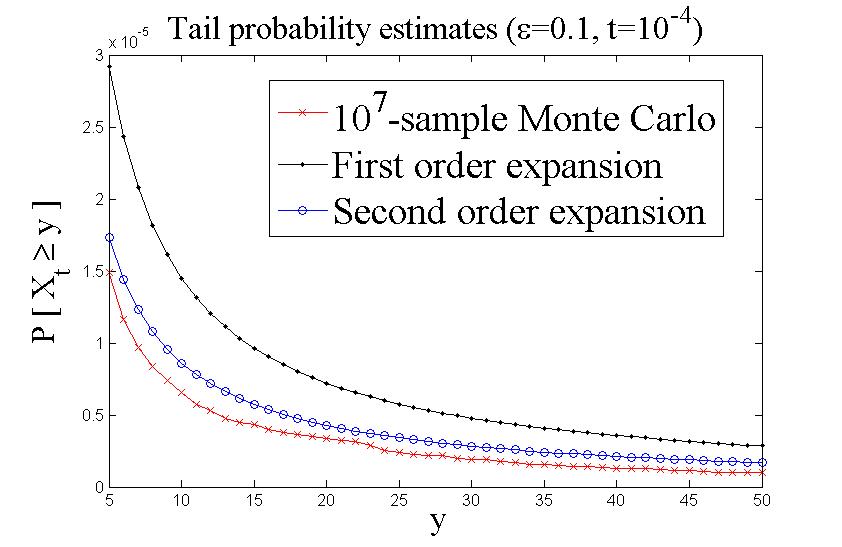}
\caption{Expansions and Monte Carlo approximation of $\mbP[X_t\geqslant y]$ with $\varepsilon=0.1$.}
 \label{ep1t4y550} 
\end{figure}

\begin{figure}[H]
\centering
\includegraphics[width=5in,height=3in]{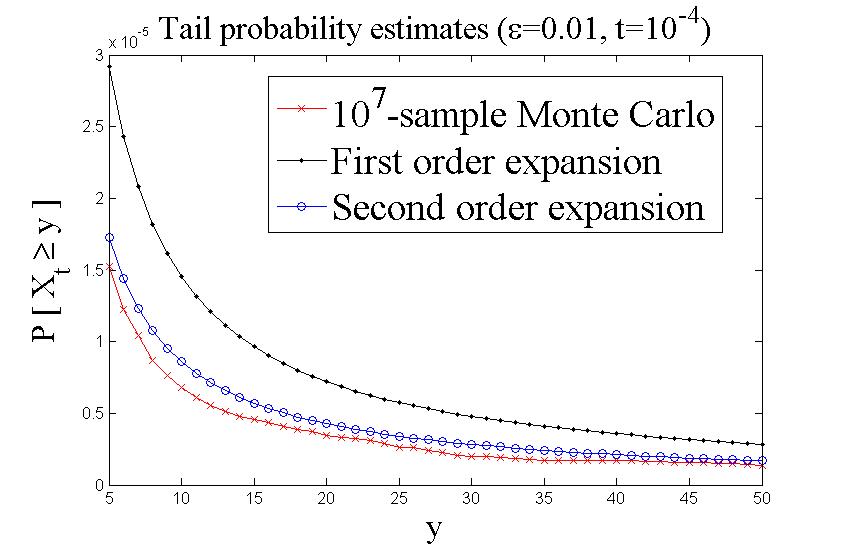}
\caption{Expansions and Monte Carlo approximation of $\mbP[X_t\geqslant y]$ with $\varepsilon=0.01$.}
 \label{ep2t4y550} 
\end{figure}

\begin{figure}[H]
\centering
\includegraphics[width=5in,height=3in]{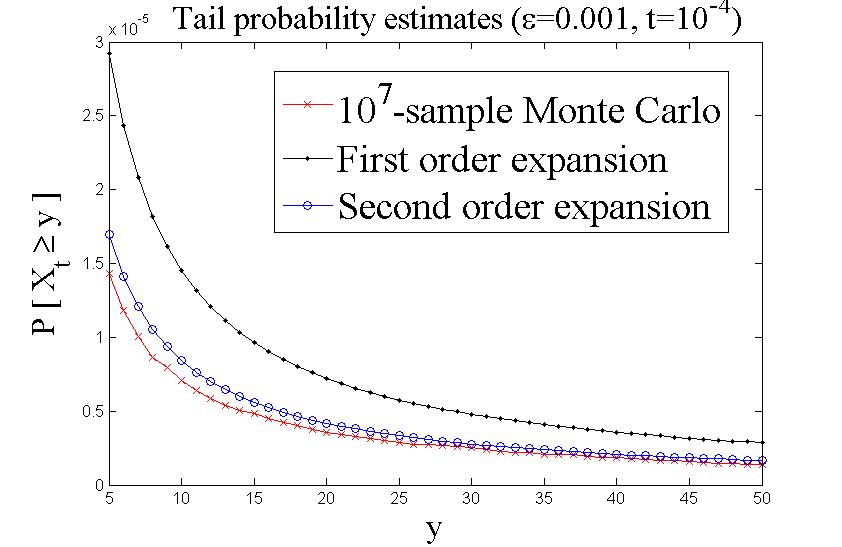}
\caption{Expansions and Monte Carlo approximation of $\mbP[X_t\geqslant y]$ with $\varepsilon=0.001$.}
 \label{ep3t4y550} 
\end{figure}

Figures \ref{ep1t4y550}-Figure \ref{ep3t4y550} show the tail probability $\mbP[X_t\geqslant{}y]$ estimated by the first and the second order expansions and the Monte Carlo approximation, for $\varepsilon=0.1,0.01,0.001$. We see that the second order expansion is indeed a better estimate than the first order expansion.
We also observe that as $\varepsilon\to0$, the second order expansion approaches to the Monte Carlo approximation and, furthermore, that the latter exhibit little variation for different $\varepsilon$, which in some sense justifies the diffusion approximation of the small-jump component of $X$.

\appendix

\section{Proofs}

\subsection{Proposition \ref{deltaproperty}}\label{pfdeltaproperty}

\noindent
\textbf{Proof:}
The first assertion holds true since
 $\lim_{w\to0^{\pm}}\psi(w)=\mp\infty$ and $\lim_{w\to\pm\infty}\bar{\psi}^{-1}(x,w)=0$ and, thus, 
 \[
 	\lim_{w\to0^{\pm}}\delta(x,w)=\gamma\left(x,\lim_{w\to{}0^{\pm}}\bar{\psi}^{-1}\left(x,\psi(w)\right)\right)=\gamma(x,0)=0.
\]
This, shows that the function $\delta$ can continuously extended on $\mbR\times(-\varepsilon,\varepsilon)$ by defining $\delta(x,0):=0$.
Now, we proceed to show the first two assertions in (\ref{1+partialxdelta}) for any $w\in(0,\varepsilon)$.
We can similarly show the case $w\in(-\varepsilon,0)$. Throughout, we use the notation $F(x,w)=\bar{\psi}^{-1}\left(x,\psi(w)\right)$ and recall that $\delta(x,w)=\gamma(x,F(x,w))$. 
 For future reference, let us also note that, for any $w\in D_{\varepsilon}$,
 \begin{equation}\label{UIFw}
 	|F(x,w)|\leq{}|w|\quad\text{and}\quad \bar\phi_{\varepsilon}(w)\leqslant\bar\phi_{\varepsilon}(F(x,w)),
\end{equation}
because
\begin{align*} \int_w^{\infty}\bar{\phi}_{\varepsilon}(r)h(r)dr\geqslant\int_{w}^{\infty}\bar{\phi}_{\varepsilon}(r)\nu(x,r)dr\quad (w>0),\quad
	\int_{-\infty}^{w}\bar{\phi}_{\varepsilon}(r)h(r)dr\geqslant\int_{-\infty}^{w}\bar{\phi}_{\varepsilon}(r)\nu(x,r)dr\quad (w<0),
\end{align*}
and the function $\bar\phi_{\varepsilon}(w)$ is chosen to be decreasing as $|w|$ increases. 

\noindent
We start to prove (\ref{1+partialxdelta}-i) for $w\in(0,\varepsilon)$. To this end, 
note that $\partial_{w}\delta(x,w)=(\partial_{2}\gamma)(x,F(x,w))\partial_{w} F(x,w)$ and, by Condition (\ref{conditiongamma}-\ref{gamma0=0}), it suffices to show that $|\partial_{w} F(x,w)|$ is bounded. By the definitions in Eq.~(\ref{DfnPhiBarPhi}) and recalling that $\nu(x,w)=\bar\nu(x,w)h(w)$ and $h(w)=g(w)|w|^{-\alpha-1}$, 
\begin{align}
  \label{partialwF0}
  \partial_wF(x,w)&=\frac{\psi'(w)}{(\partial_{2}\bar\psi)(x,F(x,w))}
  =\frac{g(w)}{\bar\nu(x,F(x,w))g(F(x,w))}\left(\frac{F(x,w)}{w}\right)^{\alpha+1}\frac{\bar\phi_{\varepsilon}(w)}{\bar\phi_{\varepsilon}(F(x,w))},
\end{align}
for  any $w\in(0,\varepsilon)$. Note that, due to (\ref{UIFw})  
 and Conditions (\ref{conditionnu}-\ref{Smothnuandbarnu}) and (\ref{conditionh}), for $\varepsilon>0$ small enough, $\inf_{x}\bar{\nu}(x,F(x,w))>\eta$ and $\inf_{x}g(F(x,w))>\eta$, for any $0<w<\varepsilon$. Using the latter and again (\ref{UIFw}), we conclude that
\[
	\sup_{0<w<\varepsilon,x\in\mbR}\left|\partial_wF(x,w)\right|<\infty,
\]
which, as discussed above, implies (\ref{1+partialxdelta}-i).

\noindent We proceed to show (\ref{1+partialxdelta}-ii), for $i=0,1,2$. The case of $i=0$ follows from  (\ref{1+partialxdelta}-i) by the mean value theorem and the fact that $\delta(x,0)=0$.
For $i=1$, note that $\partial_{x}\delta(x,w)=(\partial_{1}\gamma)(x,F(x,w))+(\partial_{2}\gamma)(x,F(x,w))\partial_{x} F(x,w)$. By Condition (\ref{conditiongamma}-\ref{gamma0=0}) and (\ref{UIFw}), 
\[
	|(\partial_{1}\gamma)(x,F(x,w))|\leq \left(\sup_{z,r}|\partial_{2}\partial_{1}\gamma(z,r)|\right)|F(x,w)|\leq{} K |w|,
\]
for some constant $K$. Similarly, since $\sup_{z,r}|\partial_{2}\gamma(z,r)|<\infty$, it suffices to show that $|\partial_{x} F(x,w)/w|$ is bounded. To this end, note that
\begin{equation}\label{partialxF0}
	\partial_{x}F(x,w)=-\frac{(\partial_{1}\bar\psi)(x,F(x,w))}{(\partial_{2}\bar\psi)(x,F(x,w))} =\frac{F(x,w)^{\alpha+1}\int_{F(x,w)}^{\infty}\partial_{1}\bar\nu(x,r)g(r)r^{-\alpha-1}\bar{\phi}_{\varepsilon}(r)dr}{\bar{\nu}(x,F(x,w))g(F(x,w))\bar\phi_{\varepsilon}(F(x,w))}.
\end{equation}
Since $\partial_{1}\bar\nu(x,r)$ and $g(r)$ are bounded for $r$ in a small neighborhood of the origin and, again, $\bar{\phi}_{\varepsilon}(w)$ is decreasing in $|w|$, for some constants $K,K'<\infty$ and any $0<w<\varepsilon$,
\begin{equation*}
	 \left|\frac{\partial_xF(x,w)}{w}\right|\leq K \frac{F(x,w)^{\alpha+1}\int_{F(x,w)}^{\infty}r^{-\alpha-1}\bar{\phi}_{\varepsilon}(r)dr}{w\bar{\nu}(x,F(x,w))g(F(x,w))\bar\phi_{\varepsilon}(F(x,w))}\leq 
	 K' \frac{F(x,w)^{\alpha+1}\bar{\phi}_{\varepsilon}(F(x,w))F(x,w)^{-\alpha}}{\bar{\nu}(x,F(x,w))g(F(x,w))\bar\phi_{\varepsilon}(F(x,w))w},
\end{equation*}
and one can use the same arguments as those in (\ref{partialwF0}) to show that $|\partial_{x} F(x,w)/w|$ is bounded on $\mbR\times(0,\varepsilon)$ and, hence, to conclude the validity of (\ref{1+partialxdelta}-i) for $i=1$. For future reference, note that the previous arguments also show that $|\partial_{x} F(x,w)/F(x,w)|$ is bounded on $\mbR\times(0,\varepsilon)$.

\noindent 
For $i=2$, note $\partial_{1}^{2}\delta(x,w)$ can be decomposed as
\[
	\left.\left[(\partial_{1}^{2}\gamma)(x,r)+
	2(\partial_{2}\partial_{1}\gamma)(x,r)\partial_{x} F(x,w)+(\partial_{2}^{2}\gamma)(x,r)\left(\partial_{x} F(x,w)\right)^{2}+(\partial_{2}\gamma)(x,r)\partial^{2}_{x} F(x,w)\right]\right|_{r=F(x,w)}.
\]
The first three terms are trivially bounded by $k |w|$, for some constant $k$, by Condition (\ref{conditiongamma}-\ref{gamma0=0}),  (\ref{UIFw}),  and the boundedness of $|\partial_{x}F(x,w)/w|$. It remains to show that $|\partial^{2}_{x}F(x,w)/w|$ is bounded. From straightforward differentiation, 
\begin{align*}
	\partial_{x}^{2}F(x,w)
	&=-\frac{(\partial_{1}^{2}\bar\psi)(x,F(x,w))}{(\partial_{2}\bar\psi)(x,F(x,w))}-2\frac{(\partial_{2}\partial_{1}\bar\psi)(x,F(x,w))\partial_{x}F(x,w)}{(\partial_{2}\bar\psi)(x,F(x,w))}\\
	&\quad +\frac{(\partial_{1}\bar\psi)(x,F(x,w))(\partial_{2}^{2}\bar\psi)(x,F(x,w))\partial_{x}F(x,w)}{\left((\partial_{2}\bar\psi)(x,F(x,w))\right)^{2}}.
\end{align*}
We denote each of the three terms on the right hand side above $D_{1},D_{2},D_{3}$, respectively, and analyze them separately. For the first term, we have:
\[
	\frac{D_{1}(x,w)}{w}=\frac{F(x,w)^{\alpha+1}\int_{F(x,w)}^{\infty}(\partial_{1}^{2}\bar\nu)(x,r)g(r)r^{-\alpha-1}\bar{\phi}_{\varepsilon}(r)dr}{\bar{\nu}(x,F(x,w))g(F(x,w))\bar\phi_{\varepsilon}(F(x,w))},
\]
which can clearly be proved to be bounded using analogous steps to those used for $\partial_{x}F(x,w)$ after Eq.~(\ref{partialxF0}). For the second term, 
\[
	\frac{D_{2}(x,w)}{w}=2\frac{(\partial_{1}\bar\nu)(x,F(x,w))}{\bar{\nu}(x,F(x,w))}\left(\frac{\partial_{x}F(x,w)}{w}\right),
\]
which is also clearly bounded since $|\partial_{x}F(x,w)/w|$ is bounded. For the third term,
\[
	\frac{D_{3}(x,w)}{w}=\frac{(\partial_{x}F(x,w))^{2}}{w}\frac{F(x,w)^{\alpha+1}\left.\partial_{r}(\bar\nu(x,r)g(r)\bar\phi_{\varepsilon}(r)r^{-\alpha-1})\right|_{r=F(x,w)}}{\bar{\nu}(x,F(x,w))g(F(x,w))\bar\phi_{\varepsilon}(F(x,w))}.
\]
The above derivatives with respect to $r$ generate the following simplified terms:
\begin{align*}
		&\frac{(\partial_{x}F(x,w))^{2}}{w}\frac{F(x,w)(\partial_{2}\bar\nu)(x,F(x,w))}{\bar{\nu}(x,F(x,w))},\quad
		\frac{(\partial_{x}F(x,w))^{2}}{wF(x,w)}\frac{F(x,w)g'(F(x,w))}{g(F(x,w))}\\
		&\frac{(\partial_{x}F(x,w))^{2}}{w}\frac{\bar\phi'_{\varepsilon}(F(x,w))}{\bar\phi_{\varepsilon}(F(x,w))},\quad 
		-(\alpha+1)\frac{\partial_{x}F(x,w)}{w}.
\end{align*}
All the above terms are clearly bounded in light of the Conditions (\ref{conditionnu}-\ref{Smothnuandbarnu}) and (\ref{conditionh}) and the facts that $|\partial_{x}F(x,w)/w|$ and  $|\partial_{x}F(x,w)/F(x,w)|$ are bounded as already proved above.

\noindent
Finally, to show (\ref{1+partialxdelta}-iii), note that, there exists an $\varepsilon>0$ small enough such that for any $|w|<\varepsilon$,
\begin{equation*}
|1+\partial_x\delta(x,w)|=|1+\partial_1\gamma(x,F(x,w))+\partial_2\gamma(x,F(x,w))\partial_xF(x,w)|>\eta,
\end{equation*}
for some $\eta>0$, due to the Condition (\ref{conditiongamma}) in Section \ref{setup} and the fact that $|\partial_{x}F(x,w)|<k|w|$, for some constant $k$.
\hfill$\blacksquare$

\subsection{Lemma \ref{counterpartA.1}}\label{pfcounterpartA.1}
\textbf{Proof:} We want to obtain the second order expansion in $t$ of the function
\begin{align}\label{DefHtz0}
  H(t,z,q)&:=\mbE\left[\bar\nu(z,J)\left.\mbP\left[Y_t(\varepsilon,\emptyset,v)\geqslant{}q\right]\right|_{v=z+\gamma(z,J)}\right],
\end{align}
where $Y(\varepsilon,\emptyset,x)$ is the solution of (\ref{Y0bigjump}), which has the same distribution law as $X^{\theta}(\varepsilon,\emptyset,x)$.
Let us start by representing (\ref{DefHtz0}) in terms of the density $\widetilde{\Gamma}(\cdot;z)$ of $z+\gamma(z,J)$, the inverse of the mapping $r\to \tilde\gamma(z,r):=z+\gamma(z,r)$, and the inverse of the {diffeomorphism} $\Phi_{t}:\eta\to Y_t(\varepsilon,\emptyset,\eta)$:
\begin{equation}\label{NdDefH}
H(t,z,q):=\mbE\left(\int \widetilde{\Gamma}(r;z)\bar\nu(z,\tilde\gamma^{-1}(z,r)){\bf 1}_{\{Y_t(\varepsilon,\emptyset,r)\geqslant{}q\}}dr\right)=\mbE\left(\int_{\Phi_{t}^{-1}(q)}^{\infty}\widetilde{\Gamma}(r;z)\bar\nu(z,\tilde\gamma^{-1}(z,r))dr\right).
\end{equation}
Let $\breve{Y}_{t}(q):=\Phi_{t}^{-1}(q)$ and recall from the proof of Lemma C.1 in \cite{Cheng} that $\breve{Y}_{t}(q)$ is a solution to the SDE
\begin{align*}
  \breve{Y}_{t}(\eta)= \eta-\int_0^t {b_{\varepsilon}(\breve{Y}_v(\eta))}dv+\int_0^t\sigma'(\breve{Y}_v(\eta))\sigma(\breve{Y}_v(\eta))dv+\int_0^t\sigma(\breve{Y}_v(\eta))d\overline{W}^T_v+{\sum^c_{0<v\leq t}}\breve{\delta}(\breve{Y}_{v-}(\eta),\Delta \bar{Z^*}^{T}_{v}),
\end{align*}
where $\breve{\delta}(u,\zeta):=\bar{\delta}(u,-\zeta)-u$ with $\bar\delta(u,\zeta)$ being the inverse of the mapping $z\to u:=z+\delta(z,\zeta)$. Above, $\{\bar{W}^{T}_{t}\}_{0\leqslant{}t\leqslant{}T}$ and $\{\bar{Z^*}^{T}_{t}\}_{0\leqslant{}t\leqslant{}T}$ are the time-reversal processes of $W^{*}$ and the L\'evy process $Z^{*}_{t}:=\int_{0}^{t}\int r\bar\mu^{*}(ds,dr)$ (see, e.g., \cite[Section VI.4]{Protter} for information about time-reversibility). In particular, $W^{T}_{\cdot}$ is a Wiener process, while $\bar{Z^*}^{T}_{\cdot}$ is an independent L\'evy process with the same law of $-Z^{*}$ (cf. Theorems VI.20 \& VI.21 in \cite{Protter}). 

\noindent
Next, we write (\ref{NdDefH}) in the form:
\[
	H(t,z,q)=\mbE\left[\widetilde{H}_{z}\left(\breve{Y}_{t}(q)\right)\right],
\]
where $\widetilde{H}_{z}(q):=\int_{q}^{\infty}\widetilde{\Gamma}(r;z)\bar\nu(z,\tilde\gamma^{-1}(z,r))dr$, and aim at applying the Dynkin's formula (\ref{Dynkin2}) to deduce:
\begin{equation}\label{DnyFrY}	
H(t,z,q)=\mbE\left[\widetilde{H}_{z}\left(\breve{Y}_{t}(q)\right)\right]
=\widetilde{H}_{z}\left(q\right) +t(\breve{L}_{\varepsilon}\widetilde{H}_{z})\left(q\right)
+t^{2}\int_{0}^{1}(1-\alpha)\mbE\left[(\breve{L}_{\varepsilon}^2\widetilde{H}_{z})\left(\breve{Y}_{\alpha{}t}(q)\right)\right]d\alpha,
\end{equation}
where $\breve{L}_{\varepsilon}$ is the infinitesimal generator of $\{\breve{Y}_t(\eta)\}_{t\geqslant0}$. We now proceed to justify (\ref{DnyFrY}) and the desired boundedness conditions (\ref{boundremainder}). To this end, it is easy to see that the following two conditions are sufficient:
\begin{equation}\label{NDCFDF2}
	{\rm (i)}\;\widetilde{H}_{z}\in C^{4}_{b}(\mbR),\qquad 
	{\rm (ii)}\;\breve{L}_{\varepsilon}f\in C^{2}_{b}(\mbR),
	\quad 
\text{for}\quad f\in C^{4}_{b}(\mbR).
\end{equation}
The condition (\ref{NDCFDF2}-i) follows directly from Condition (\ref{conditionnu}) and Lemma \ref{gammainfinity}. For the other condition in (\ref{NDCFDF2}), let us note that (see \cite[Eq.~(3.5)]{Cheng})
\begin{align}\label{generatorsmalla}
\breve{L}_{\varepsilon}f(y)&:=\breve{b}_{\varepsilon}(y)f'(y)+\frac{\sigma^2(y)}{2}f''(y)+\mathcal{I}_{\varepsilon}f(y),\quad\text{where}\\
\breve{\mathcal{I}}_{\varepsilon}f(y)
&=\int\int_{0}^{1}f''(y+\breve\delta(y,r)\beta)(1-\beta)d\beta\left(\breve{\delta}(y,r)\right)^{2}\bar{h}_{\varepsilon}(-r)dr\label{Irepresentation2a}.
\end{align}
and $\breve{b}_{\varepsilon}(y):=-b_{\varepsilon}(y)+\sigma'(y)\sigma(y)$. Clearly, the first two terms in $\breve{L}_{\varepsilon}f$ belong to $C^{2}_{b}(\mbR)$ when $f\in C^{4}_{b}$, provided that $\sigma^{2},\breve{b}_{\varepsilon}\in C^{2}_{b}$, which holds true in light of Conditions (\ref{conditionnu}), (\ref{Cinfinity}-\ref{bsigma}), and (\ref{conditiongamma}-\ref{gamma0=0}). The boundedness of $\breve{\mathcal{I}}_{\varepsilon}f$ would hold true provided that $f\in C^{2}_{b}$ and 
\begin{equation}\label{Cnd0delta}
	\sup_{u}|\breve{\delta}(u,r)|\leq{}k |r|, \quad r\in (-\varepsilon,\varepsilon),
\end{equation}
for some constant $k$. Since $\breve{\delta}(u,0)=0$, it suffices to show that $\sup_{u\in\mbR,r\in D_{\varepsilon}}|\partial_{r}\breve{\delta}(u,r)|<\infty$, where $D_{\varepsilon}:=(-\varepsilon,0)\cup(0,\varepsilon)$. Note that 
\[
	\partial_{r}\breve{\delta}(u,r)=-(\partial_{2}\bar\delta)(u,-r)=\frac{(\partial_{2}\delta)(\bar\delta(u,-r),-r)}{1+(\partial_{1}\delta)(\bar\delta(u,-r),-r)}, \quad r\in D_{\varepsilon},
\]
which is bounded on $\mbR\times D_{\varepsilon}$ by Proposition \ref{deltaproperty}. The formal differentiation of $\breve{\mathcal{I}}_{\varepsilon}f(y)$ yields
\begin{align}\label{DerGenYb}
	\partial_{y}\breve{\mathcal{I}}_{\varepsilon}f(y)&=\int\int_{0}^{1}f^{(3)}(y+\breve\delta(y,r)\beta)(1-\beta)\left(1+(\partial_{1}\breve\delta)(y,r)\beta\right)d\beta\left(\breve{\delta}(y,r)\right)^{2}\bar{h}_{\varepsilon}(-r)dr\\
	&\quad +2\int\int_{0}^{1}f''(y+\breve\delta(y,r)\beta)(1-\beta)d\beta\breve{\delta}(y,r)(\partial_{1}\breve{\delta})(y,r)\bar{h}_{\varepsilon}(-r)dr,\nonumber
\end{align}
and it is clear that for the derivative to be well-defined and bounded it suffices that $f\in C_{b}^{3}(\mbR)$, (\ref{Cnd0delta}), and 
\begin{equation}\label{Cnd1delta}
	\sup_{u}|(\partial_{1}\breve{\delta})(u,r)|\leq{}k |r|, \quad r\in (-\varepsilon,\varepsilon),
\end{equation}
Note that 
\[
	(\partial_{1}\breve{\delta})(u,r)=(\partial_{1}\bar\delta)(u,-r)-1=\frac{(\partial_{1}\delta)(\bar\delta(u,-r),-r)}{1+(\partial_{1}\delta)(\bar\delta(u,-r),-r)}, \quad r\in D_{\varepsilon},
\]
which can be bounded by $k|r|$ on $\mbR\times D_{\varepsilon}$ by Proposition \ref{deltaproperty}. 

\noindent 
Formally differentiating (\ref{DerGenYb}), generates the following terms:
\begin{align}\label{TwiceDerGenYb}
	\partial_{y}\breve{\mathcal{I}}_{\varepsilon}f(y)&=\int\int_{0}^{1}f^{(4)}(y+\breve\delta(y,r)\beta)(1-\beta)\left(1+(\partial_{1}\breve\delta)(y,r)\beta\right)^{2}d\beta\left(\breve{\delta}(y,r)\right)^{2}\bar{h}_{\varepsilon}(-r)dr\\
	&\quad+\int\int_{0}^{1}f^{(3)}(y+\breve\delta(y,r)\beta)(1-\beta)(\partial_{1}^{2}\breve\delta)(y,r)\beta d\beta\left(\breve{\delta}(y,r)\right)^{2}\bar{h}_{\varepsilon}(-r)dr\nonumber
	\\
	&\quad +4\int\int_{0}^{1}f^{(3)}(y+\breve\delta(y,r)\beta)(1-\beta)\left(1+(\partial_{1}\breve\delta)(y,r)\beta\right)d\beta\breve{\delta}(y,r)(\partial_{1}\breve{\delta})(y,r)\bar{h}_{\varepsilon}(-r)dr\nonumber
	\\
		&\quad +2\int\int_{0}^{1}f''(y+\breve\delta(y,r)\beta)(1-\beta)d\beta\left((\partial_{1}\breve{\delta})(y,r)\right)^{2}\bar{h}_{\varepsilon}(-r)dr,\nonumber\\
		&\quad +2\int\int_{0}^{1}f''(y+\breve\delta(y,r)\beta)(1-\beta)d\beta \breve{\delta}(y,r)(\partial_{1}^{2}\breve{\delta})(y,r)\bar{h}_{\varepsilon}(-r)dr.\nonumber
\end{align}
The previous expression shows that in order for the second derivative to be well-defined and bounded it suffices that both conditions (\ref{Cnd0delta}) and (\ref{Cnd1delta}) are satisfied as well as 
\begin{equation}\label{Cnd2delta}
	\sup_{u}|(\partial_{1}^{2}\breve{\delta})(u,r)|\leq{}k |r|, \quad r\in (-\varepsilon,\varepsilon).
\end{equation}
The latter condition again follows from Proposition \ref{deltaproperty} since, for $r\in D_{\varepsilon}$,
\[
	(\partial_{1}^{2}\breve{\delta})(u,r)=-\frac{(\partial_{1}^{2}\delta)(\bar\delta(u,-r),-r)}{\left(1+(\partial_{1}\delta)(\bar\delta(u,-r),-r)\right)^{3}}.
\]

\noindent
Some further algebra shows that the terms $\widetilde{H}_{z}\left(q\right)$ and $(\breve{L}_{\varepsilon}\widetilde{H}_{z})\left(q\right)$ in (\ref{DnyFrY}) coincide with expressions of $H_{0}(z;q)$ and $H_{1}(z,q)$ given in  (\ref{Exp1J}). From (\ref{NDCFDF2}) and the fact that the derivatives $\widetilde{H}^{(i)}_{z}(q)$, $i=0,\dots,4$, are bounded both on $z$ and $q$, it is clear that $\sup_{z,q}|(\breve{L}_{\varepsilon}\widetilde{H}_{z})\left(q\right)|<\infty$ and $\sup_{z,q}|(\breve{L}_{\varepsilon}^{2}\widetilde{H}_{z})\left(q\right)|<\infty$, which in turn imply the conditions in Eq.~(\ref{boundremainder}).
 \hfill$\blacksquare$

\subsection{Lemma \ref{Y0diffeomorphism}}\label{pfY0diffeomorphism}
\noindent\textbf{Proof:} Throughout the proof, we write $Y_{t}^{x}:=Y_t(\varepsilon,\emptyset,z)$. Let us start by noting that the SDE (\ref{Y0bigjump}) satisfies the Hypotheses 5-9 in \cite{Malliavin}, with $E=(-\varepsilon,\varepsilon)$, $G(dr)=\bar{\phi}_{\varepsilon}(r)h(r)dr$, and $\eta(r)=r$, since obviously $\int_{-\varepsilon}^{\varepsilon}|\eta(r)|^ph(r)dr<\infty$, for any $p\geqslant2$, and, by  Theorem \ref{deltaproperty},  both
\begin{align*}
  \left|\frac{\delta(z,r)}{r}\right|
  \quad\text{\rm and}\quad
  \left|\frac{\partial_z\delta(z,r)}{r}\right|
\end{align*}
are bounded.
Then, Theorem 5-10 in \cite{Malliavin} guarantees the existence of a unique solution of (\ref{Y0bigjump}).
In turn the latter shows the existence of a unique weak solution for the SDE (\ref{X0bigjump}). Using an interlacing technique similar to that used in Theorem 6.2.9 of \cite{Applebaum}, one can proceed to show the existence of a unique weak solution to (\ref{Xtheta}), which in turn implies the existence of a unique weak solution for the SDE (\ref{SDEpq}).

\noindent
To conclude that $x\to{}Y_{t}^{x}$ is a diffeomorphism, let us frame the family of solutions $\{Y_t^{z}\}_{z\in\mbR}$ of the SDE (\ref{Y0bigjump}) in the form (5-22) of \cite{Malliavin}, indexed by the initial state $z$ in a bounded neighborhood $U$ of $x\in\mbR$, with $H^{z}=z$, $A^{z}=b_{\varepsilon}$, $B^{z}=\sigma$ and $C^{z}=\delta$. First, note that the SDE (\ref{Y0bigjump}) satisfies the conditions \rm(i), (ii), and (iv) of \cite[Hypothesis 5-23]{Malliavin} since the coefficients $b_{\varepsilon}(y)$, $\sigma(y)$ and $\delta(y,r)/r$ are twice differentiable in $y$ and their respective partial derivatives with respect to $y$ are bounded in light of Assumptions (\ref{conditionnu})-(\ref{conditiongamma}) and Proposition \ref{deltaproperty}. The assumption {\rm(iii)} of \cite[Hypothesis 5-23]{Malliavin}, and {\rm(i)-(ii)} in \cite[Theorem 5-24]{Malliavin} are trivially satisfied since $A^z=b_{\varepsilon},B^z=\sigma,C^z=\delta$ are deterministic functions independent of $z$. Therefore, by Theorem  5.24 in \cite{Malliavin}, the mapping $\Phi:x\mapsto{}Y_t^{x}$ is differentiable and its derivative, $\partial_x Y_{t}^{x}$ is the unique solution of the SDE obtained by formal differentiation of (\ref{Y0bigjump}):
\begin{align}\label{partialxY}
  \partial_xY^x_t
  &=1+\int_0^s b'_{\varepsilon}\left(Y_v^x\right)\partial_xY^x_{v}dv
  +\int_0^t\sigma'\left(Y_v^x\right)\partial_xY^x_{v}dW^{*}_v
  +\int_0^s\int\partial_1\delta\left(Y_{v^{-}}^x,r\right)\partial_xY^x_{v^-}\bar{\mu}^{*}(dv,dr),
\end{align}
Note that the coefficients of (\ref{Y0bigjump}) are deterministic functions instead of stochastic processes, so the terms $\partial{}A, \partial{}B$ and $\partial{}C$ in (5-25) of \cite{Malliavin} are absent in (\ref{partialxY}) above. In particular, $\partial_xY^x_t$ is given by
\[
	 \partial_xY^x_t=\exp\left(V_{t}-\frac{1}{2}\int_0^t \sigma'\left(Y_v^x\right)^{2}dv-\int_{0}^{t}\int\ln\left(1+\partial_1\delta\left(Y_{v^{-}}^x,r\right)\right){\mu}^{*}(dv,dr)\right),
 \]
 where $V_{t}:=\int_0^t b'_{\varepsilon}\left(Y_v^x\right)dv
  +\int_0^t\sigma'\left(Y_v^x\right)dW^{*}_v
  +\int_0^t\int\partial_1\delta\left(Y_{v^{-}}^x,r\right)\bar{\mu}^{*}(dv,dr)$. Due to (\ref{1+partialxdelta}), it is clear that, a.s.,  $\partial_xY^x_t\neq{} 0$, for all $t$. Hence, the mapping $\Phi$ admits a differentiable inverse by the implicit function theorem, and, finally, $\Phi$ is a diffeomorphism on $\mbR$.
\hfill$\blacksquare$

\subsection{Theorem \ref{tailexpansion}}\label{pftailexpansion}

\textbf{Proof: }Throughout the proof, we only assume that $\nu(x,r)\leqslant{}h(r)$ and the function $\bar{\nu}(x,r):=\nu(x,r)/h(r)$ admits an extension on $\mbR\times\mbR$ that is $C_b^4$ in $x$, which is weaker than the technical assumption (\ref{conditionnu}) in Section \ref{setup}.

\noindent\textbf{The case $N_t^{\varepsilon}=0$.}

\noindent Recalling that we denote the law of a process $V$ (respectively, the conditional law of $V$ given an event $B$) by $\mathcal{L}(V)$ (resp. $\mathcal{L}(V|B)$), we have
\begin{align*}
\mathcal{L}\left(\left.\left\{X^{\theta}_s(\varepsilon,x)\right\}_{0\leqslant{}s\leqslant{}t}\right|N_t^{\varepsilon}=0\right)=\mathcal{L}\left(\left\{X^{\theta}_s(\varepsilon,\emptyset,x)\right\}_{0\leqslant{}s\leqslant{}t}\right).
\end{align*}
Consequently, the inequality (\ref{ttoN}) implies that
\begin{align*}
\mbP\left[\left.X^{\theta}_{t}(\varepsilon,x)\geqslant{}x+y\right|N_{t}^{\varepsilon}=0\right]\leqslant{}Ct^N,
\end{align*}
for any $y>0$, where $N>0$ can be made arbitrarily large by taking $\varepsilon>0$ small enough.

\noindent\textbf{The case $N_t^{\varepsilon}=1$.}

\noindent Conditioning on the time of the jump, we have
\begin{align*}
\mbP\left[\left.X^{\theta}_{t}(\varepsilon,x)\geqslant{}x+y\right|N_{t}^{\varepsilon}=1\right]=\frac{1}{t}\int_0^t\mbP\left[X^{\theta}_t(\varepsilon,\{s\},x)\geqslant{}x+y\right]ds.
\end{align*}
Let $J$ and $U$ denote, respectively, a random variable with density function $\breve{h}_{\varepsilon}(r)=h_{\varepsilon}(r)/\lambda_{\varepsilon}$ and an independent random variable with standard uniform distribution on $(0,1)$.  Then, conditioning on $\mathcal{F}_{s^-}$, the integrand above can be written as follows:
\begin{align*}
\mbP\left[X^{\theta}_t(\varepsilon,\{s\},x)\geqslant{}x+y\right]
&=\mbE\left[\mbP\left[X^{\theta}_{t-s}\left(\varepsilon,\emptyset,v\right)\geqslant{}x+y\right]_{v=\hat{X}_{s^-}(\varepsilon,\emptyset,x)+\gamma\left(\hat{X}_{s^-}(\varepsilon,\emptyset,x),J\right)\theta\left(\hat{X}_{s^-}(\varepsilon,\emptyset,x),J,U\right)}\right]\\
&=\mbE\left[\mathbf{1}_{U>\bar{\nu}\left(\hat{X}_{s^-}(\varepsilon,\emptyset,x),J\right)}\mbP\left[X^{\theta}_{t-s}\left(\varepsilon,\emptyset,v\right)\geqslant{}x+y\right]_{v=\hat{X}_{s^-}(\varepsilon,\emptyset,x)}\right]\\
&\quad+\mbE\left[\mathbf{1}_{U<\bar{\nu}\left(\hat{X}_{s^-}(\varepsilon,\emptyset,x),J\right)}\mbP\left[X^{\theta}_{t-s}\left(\varepsilon,\emptyset,v\right)\geqslant{}x+y\right]_{v=\hat{X}_{s^-}(\varepsilon,\emptyset,x)+\gamma\left(\hat{X}_{s^-}(\varepsilon,\emptyset,x),J\right)}\right],\\
&=\mbE\left[\left(1-\bar{\nu}\left(\hat{X}_{s^-}(\varepsilon,\emptyset,x),J\right)\right)\mbP\left[X^{\theta}_{t-s}\left(\varepsilon,\emptyset,v\right)\geqslant{}x+y\right]_{v=\hat{X}_{s^-}(\varepsilon,\emptyset,x)}\right]\\
&\quad+\mbE\left[\bar{\nu}\left(\hat{X}_{s^-}(\varepsilon,\emptyset,x),J\right)\mbP\left[X^{\theta}_{t-s}\left(\varepsilon,\emptyset,v\right)\geqslant{}x+y\right]_{v=\hat{X}_{s^-}(\varepsilon,\emptyset,x)+\gamma\left(\hat{X}_{s^-}(\varepsilon,\emptyset,x),J\right)}\right].
\end{align*}
where $\{\hat{X}_{s}(\varepsilon,\emptyset,x)\}_{0\leqslant{}s\leqslant{}t}$ is an independent copy of $\{X^{\theta}_{s}(\varepsilon,\emptyset,x)\}_{0\leqslant{}s\leqslant{}t}$. We denote the last two terms above by $T_1(x,y)$ and $T_2(x,y)$, respectively. By the fact that $\nu\leqslant{}h$ and from the Markov's property,
\[
T_{1}(t)\leqslant\mbE\left[\left.\mbP\left(X^{\theta}_{t-s}\left(\varepsilon,\emptyset,z\right)\geqslant{}x+y\right)\right|_{z=\hat{X}_{s^-}(\varepsilon,\emptyset,x)}\right]=\mbP[X^{\theta}_t(\varepsilon,\emptyset,x)\geqslant{}x+y],
\]
which can be made $O(t^N)$ for any $y>0$ and an arbitrary large $N>0$ in light of (\ref{ttoN}). On the other hand,
\begin{equation*}
T_2(x,y)=\mbE\left[H\left(t-s,\hat{X}_{s^-}(\varepsilon,\emptyset,x),x+y\right)\right].
\end{equation*}
where
\begin{align*}
  H(t,z,q)=\mbE\left[\bar{\nu}(z,J)\mbP\left[X_t^{\theta}(\varepsilon,\emptyset,v)\geqslant{}q\right]_{v=z+\gamma(z,J)}\right].
\end{align*}
Using Theorem \ref{counterpartA.1}, $T_2(x,y)$ can be written as
\begin{equation}\label{T2t}
T_2(x,y)=\mbE\left[H_0\left(\hat{X}_{s^-}(\varepsilon,\emptyset,x);x+y\right)+(t-s)H_1\left(\hat{X}_{s^-}(\varepsilon,\emptyset,x);x+y\right)\right]+O(t^2),
\end{equation}
where $H_0(z;q)$ and $H_1(z;q)$ are given in (\ref{Exp1J}).
By writing $H_{0}(z;q)$ as $\lambda_{\varepsilon}^{-1}\int^{\infty}_{\gamma^{-1}(z,q-z)}\bar\nu(z,r)h_{\varepsilon}(r)dr$ and recalling from Lemma \ref{gammainfinity} that $\gamma^{-1}$ is $\bar{C}^{4}_{b}$ and the regularity on $\bar\nu$ imposed by (\ref{conditionnu}), it follows that $H_0(z;q)$ is $C_b^4$ in $z$. Then we can apply the second order Dynkin's formula (\ref{Dynkin2}) to $\mbE[H_0(\hat{X}_{s}(\varepsilon,\emptyset,x);x+y)]$ to get
\begin{equation}\label{EH0}
\mbE\left[H_0\left(\hat{X}_{s}(\varepsilon,\emptyset,x);x+y\right)\right]
=H_{0,0}(x;y)+sH_{0,1}(x;y)+O(s^2),
\end{equation}
where
\begin{align}
  \label{H00}H_{0,0}(x;y)&:=H_{0}(x;x+y)
  =\frac{1}{\lambda_{\varepsilon}}\int^{\infty}_{\gamma^{-1}(x,y)}\nu_{\varepsilon}(x,r)dr,\\
 	H_{0,1}(x;y)&:=L_{\varepsilon}H_0(x;x+y)
  =b_{\varepsilon}(x)\partial_1H_0(x;x+y)
  +\frac{\sigma^2(x)}{2}\partial_1^2H_0(x;x+y)
  +\hat{H}_{0,1}(x;y),\label{H01}
\end{align}
with
\begin{align*}
  \hat{H}_{0,1}(x;y)=\int\left[H_0(x+\gamma(x,r);x+y)-H_0(x;x+y)-\partial_1H_0(x;x+y)\gamma(x,r)\right]\bar\nu_{\varepsilon}(x,r)dr.
\end{align*}
From the definition of $H_{0}(z;q):=\lambda_{\varepsilon}^{-1}\int^{\infty}_{\gamma^{-1}(z,q-z)}\nu_{\varepsilon}(z,r)dr$, one can readily check that
\begin{align}
  \label{partial1H0xy}
  \lambda_{\varepsilon}\partial_1H_0(x;x+y)
  &=-\nu_{\varepsilon}\left(x,\gamma^{-1}(x,y)\right)\left[\partial_1\gamma^{-1}(x,y)-\partial_2\gamma^{-1}(x,y)\right]
  +\int^{\infty}_{\gamma^{-1}(x,y)}\partial_1\nu_{\varepsilon}(x,r)dr,\\
  \label{partial1secondH0xy}
  \lambda_{\varepsilon}\partial_1^2H_0(x;x+y)&=-\partial_2\nu_{\varepsilon}\left(x,\gamma^{-1}(x,y)\right)\left(\partial_1\gamma^{-1}(x,y)-\partial_2\gamma^{-1}(x,y)\right)^{2}+\int^{\infty}_{\gamma^{-1}(x,y)}\partial_1^2\nu_{\varepsilon}(x,r)dr\\
  &\quad -\nu_{\varepsilon}\left(x,\gamma^{-1}(x,y)\right)\left[\partial_1^{2}\gamma^{-1}(x,y)-2\partial_{1}\partial_{2}\gamma^{-1}(x,y)+\partial_2^{2}\gamma^{-1}(x,y)\right]\nonumber\\
  \label{hatH01}\lambda_{\varepsilon}\hat{H}_{0,1}(x;y)&=\int\bigg(\int^{\infty}_{\gamma^{-1}(x+\gamma(x,r),y-\gamma(x,r))}\nu_{\varepsilon}(x+\gamma(x,r),r_1)dr_1\\
  \nonumber&\qquad\quad-\int_{\gamma^{-1}(x,y)}^{\infty}\nu_{\varepsilon}(x,r_1)dr_1
  -\lambda_{\varepsilon}\partial_1H_0(x;x+y)\gamma(x,r)\bigg)\bar{\nu}_{\varepsilon}(x,r)dr.
\end{align}
Next, in order to apply the first order Dynkin's formula (\ref{Dynkin1}) to $\mbE[H_1(\hat{X}_s(\varepsilon,\emptyset,x);x+y)]$, we need to first show that the mapping $z\mapsto{}H_1(z;x+y)$ belongs to $C_b^2$. Since $H_1(z;q)=D(z;q)+I(z;q)$ with $D(\cdot;q)\in{}C_b^2$ by the regularity of $\nu(z,r)$ in $z$ and by Lemma \ref{gammainfinity}-1, we only need to show that $I(\cdot;q)\in{}C_b^2$. To this end, let us write
\begin{equation}\label{ReprIEsy}
	I(z;q)=\int \int_{0}^{1}\partial_{2}^{2}\mathcal{C}(z,\beta r;q)(1-\beta)d\beta h(r)\bar{\phi}_{\varepsilon}(r) r^{2}dr,
\end{equation}
where
\begin{align*}
  \mathcal{C}(z,r;q):=\int^q_{\bar{\gamma}(q,r)}\nu_{\varepsilon}\left(z,\tilde{\gamma}^{-1}(z,\eta)\right)\partial_2\tilde{\gamma}^{-1}(z,\eta)\bar\nu(\eta,r)d\eta-\nu_{\varepsilon}\left(z,\tilde{\gamma}^{-1}(z,q)\right)\partial_2\tilde{\gamma}^{-1}(z,q)\gamma(q,r)\bar\nu(q,r)
\end{align*}
The representation (\ref{ReprIEsy}) follows by the Taylor's theorem and the fact that $\mathcal{C}(z,0;q)=\partial_2\mathcal{C}(z,0;q)=0$. Next, observe that $\mathcal{C}(z,r;q)$ is $C^{2}_{b}$ in $(z,r)$ since all the involved functions are $C^{2}_{b}$. Therefore, standard applications of the dominated convergence theorem implies that $I(z;q)$ is $C_{b}^{2}$ in $z$ and, furthermore,
\begin{align*} \partial^{i}_{1}I(z;q)=\int\int_{0}^{1}\partial^{i}_{1}(\partial_{2}^{2}\mathcal{C}(z,\beta{}r;q))(1-\beta)d\beta{}h(r)\bar{\phi}_{\varepsilon}(r)r^{2}dr,
\quad{}i=1,2.
\end{align*}
Then, we apply the first order Dynkin's formula (\ref{Dynkin1}) to $\mbE[H_1(\hat{X}_s(\varepsilon,\emptyset,x);x+y)]$:
\begin{align}\label{EH1}
\mbE\left[H_1\left(\hat{X}_{s^-}(\varepsilon,\emptyset,x);x+y\right)\right]
=H_{1,0}(x;y)+O(s),
\end{align}
where, using the following relationships obtained by implicit function theorem
\begin{align*}
  \tilde\gamma^{-1}(x,x+y)=\gamma^{-1}(x,y),
  \quad\partial_2\tilde{\gamma}^{-1}(x,x+y)=\partial_2\gamma^{-1}(x,y),
  \quad\partial^2_2\tilde{\gamma}^{-1}(x,x+y)=\partial_2^2\gamma^{-1}(x,y),
\end{align*}
we have
\begin{equation}\label{H10}
\begin{aligned}
H_{1,0}(x;y)&:=H_{1}(x;x+y)
=D_{1,0}(x;y)+I_{1,0}(x;y),\quad\text{with}\\	
D_{1,0}(x;y)&:=\frac{1}{\lambda_{\varepsilon}}\Big(\big[b_{\varepsilon}(x+y)-v'(x+y)\big]\nu_{\varepsilon}\left(x,\gamma^{-1}(x,y)\right)\partial_2\gamma^{-1}(x,y)\\
&\qquad\quad-v(x+y)\big[\partial_{2}\nu_{\varepsilon}\left(x,{\gamma}^{-1}(x,y)\right)\big[\partial_2{\gamma}^{-1}(x,y)\big]^{2}
+\nu_{\varepsilon}\left(x,{\gamma}^{-1}(x,y)\right)\partial^{2}_2{\gamma}^{-1}(x,y)\big]\Big),\\
I_{1,0}(x;y)&:=\frac{1}{\lambda_{\varepsilon}}\int\bigg[\int^{x+y}_{\bar{\gamma}(x+y,r)}\nu_{\varepsilon}\left(x,\tilde{\gamma}^{-1}(x,\eta)\right)\partial_2\tilde{\gamma}^{-1}(x,\eta)\bar\nu_{\varepsilon}(\eta,r)d\eta\\
  &\quad\qquad\qquad-\nu_{\varepsilon}\left(x,\gamma^{-1}(x,y)\right)\partial_2\gamma^{-1}(x,y)\gamma(x+y,r)\bar\nu_{\varepsilon}(x+y,r)\bigg]dr.
\end{aligned}
\end{equation}
Finally, substitute (\ref{EH0}) and (\ref{EH1}) into (\ref{T2t}), we have
\begin{align}\label{1JExp}
  \mbP\left[X^{\theta}_t(\varepsilon,\{s\},x)\geqslant{}x+y\right]
  =H_{0,0}(x;y)+sH_{0,1}(x;y)+(t-s)H_{1,0}(x;y)+O(t^2),
\end{align}
where $H_{0,0}, H_{0,1}$ and $H_{1,0}$ are given in (\ref{H00}), (\ref{H01}) and (\ref{H10}), respectively. Thus,
\begin{align}\label{ConN1}
\mbP\left[\left.X^{\theta}_{t}(\varepsilon,x)\geqslant{}x+y\right|N_{t}^{\varepsilon}=1\right]
=H_{0,0}(x;y)+\frac{t}{2}[H_{0,1}(x;y)+H_{1,0}(x;y)]+O(t^2).
\end{align}

\noindent\textbf{The case $N_t^{\varepsilon}=2$}

\noindent
Conditioning on the times of the jumps, we have
\begin{align}\label{N2int}
\mbP\left[\left.X^{\theta}_{t}(\varepsilon,x)\geqslant{}x+y\right|N_{t}^{\varepsilon}=2\right]=\frac{2}{t^2}\int_0^t\int_{s_1}^t\mbP\left[X^{\theta}_t(\varepsilon,\{s_1,s_2\},x)\geqslant{}x+y\right]ds_2ds_1.
\end{align}
Denote the probability inside the integral above by $A_{t}(s_{1},s_{2},x,y)$.
Below, let $J_i$ and $U_i$ ($i=1,2$) be independent copies of $J$ and $U$, respectively. Denoting an independent copy of $\{X^{\theta}_{s}(\varepsilon,\{s_1\},x)\}_{0\leqslant{}s\leqslant{}t}$ by $\{\hat{X}_{s}(\varepsilon,\{s_1\},x)\}_{0\leqslant{}s\leqslant{}t}$  and conditioning on $\mathcal{F}_{s_2^-}$, we have
\begin{align}
\nonumber{}A_{t}(s_{1},s_{2},x,y)
&=\mbE\left[\mbP\left[X^{\theta}_{t-s_2}\left(\varepsilon,\emptyset,v\right)\geqslant{}x+y\right]_{v=\hat{X}_{s_2^-}(\varepsilon,\{s_1\},x)+\gamma\left(\hat{X}_{s_2^-}(\varepsilon,\{s_1\},x),J\right)\theta\left(\hat{X}_{s_2^-}(\varepsilon,\{s_1\},x),J_2,U_2\right)}\right]\\
\label{2bigjumps}&=\mbE\left[\Psi_{s_2,t}\left(\hat{X}_{s_2}(\varepsilon,\{s_1\},x)\right)\right]
+\mbE\left[\bar{\Psi}_{s_2,t}\left(\hat{X}_{s_2}(\varepsilon,\{s_1\},x)\right)\right].
\end{align}
where
\begin{align*} \Psi_{s_2,t}(x_{1})&:=\mbE\left[1-\bar{\nu}\left(x_{1},J_{2}\right)\right]\mbP\left[X^{\theta}_{t-s_2}\left(\varepsilon,\emptyset,x_{1}\right)\geqslant{}x+y\right]
=\mbE\big[\big(1-\bar{\nu}(x_{1},J_{2})\big)\mathbf{1}_{X^{\theta}_{t-s_2}\left(\varepsilon,\emptyset,x_{1}\right)\geqslant{}x+y}\big],\\ \bar\Psi_{s_2,t}(x_{1})&:=\mbE\left[\bar{\nu}\left(x_{1},J_{2}\right)\left.\mbP\left[X^{\theta}_{t-s_2}\left(\varepsilon,\emptyset,v\right)\geqslant{}x+y\right]\right|_{v=x_{1}+\gamma(x_{1},J_{2})}\right],
\end{align*}
We denote the last two terms on the right-hand side of (\ref{2bigjumps}) by $T_3(x,y)$ and $T_4(x,y)$, respectively.
Conditioning on $\mathcal{F}_{s_1^-}$, we can write $T_3$ as
\begin{align}
\nonumber
T_3&=\mbE\left[\left.\mbE\left[\Psi_{s_2,t}\left(\hat{X}_{s_2-s_1}(\varepsilon,\emptyset,x_{2})\right)\right]\right|_{x_{2}=\tilde{X}_{s_1}(\varepsilon,\emptyset,x)+\gamma(\tilde{X}_{s_1}(\varepsilon,\emptyset,x),J_1)\theta(\tilde{X}_{s_1}(\varepsilon,\emptyset,x),J_1,U_1)}\right]\\ \nonumber&=\mbE\left[\left(1-\bar\nu\left(\tilde{X}_{s_1}(\varepsilon,\emptyset,x),J_{1}\right)\right)\left.\mbE\left[\Psi_{s_2,t}\left(\hat{X}_{s_2-s_1}(\varepsilon,\emptyset,x_{2})\right)\right]\right|_{x_{2}=\tilde{X}_{s_1}(\varepsilon,\emptyset,x)}\right]\\ &\quad+\mbE\left[\bar\nu\left(\tilde{X}_{s_1}(\varepsilon,\emptyset,x),J_{1}\right)\left.\mbE\left[\Psi_{s_2,t}\left(\hat{X}_{s_2-s_1}(\varepsilon,\emptyset,x_{2})\right)\right]\right|_{x_{2}=\tilde{X}_{s_1}(\varepsilon,\emptyset,x)+\gamma(\tilde{X}_{s_1}(\varepsilon,\emptyset,x),J_1)}\right] \nonumber\\
	&=:T_{3,1}+T_{3,2},\label{T3inner2}
\end{align}
where $\{\tilde{X}_s(\varepsilon,\emptyset,x)\}_{0\leqslant{}s\leqslant{}t}$ and $\{\hat{X}_s(\varepsilon,\emptyset,x)\}_{0\leqslant{}s\leqslant{}t}$ are independent copies of $\{X^{\theta}_s(\varepsilon,\emptyset,x)\}_{0\leqslant{}s\leqslant{}t}$. Using that $0\leqslant{}\bar\nu\leqslant{}1$ and the Markov's property, it is clear that
\begin{align}
\nonumber|T_{3,1}|&\leqslant\mbE\left[\left.\mbE\left[\left|\Psi_{s_2,t}\left(\hat{X}_{s_2-s_1}(\varepsilon,\emptyset,x_{2})\right)\right|\right]\right|_{x_{2}=\tilde{X}_{s_1}(\varepsilon,\emptyset,x)}\right]\\
\nonumber&=\mbE\left[\left|\Psi\left(\hat{X}_{s_2}(\varepsilon,\emptyset,x)\right)\right|\right]\\
\nonumber&\leqslant\mbE\left[\left.\mbP\left[X^{\theta}_{t-s_2}\left(\varepsilon,\emptyset,x_{1}\right)\geqslant{}x+y\right]\right|_{x_{1}=\hat{X}_{s_2}(\varepsilon,\emptyset,x)}\right]\\
&=\mbP\left[X^{\theta}_{t}\left(\varepsilon,\emptyset,x\right)\geqslant{}x+y\right]=O(t),\label{BndT31}
\end{align}
as $t\to0$, by taking $\varepsilon>0$ small enough. To deal with the second term in (\ref{T3inner2}), let us define
\begin{align*} \widetilde{\Psi}_{s_{2},t}(x_{1}):=\mbE\left[\left(1-\bar{\nu}\left(X^{\theta}_{t-s_2}\left(\varepsilon,\emptyset,x_{1}\right),J_{2}\right)\right)\mathbf{1}_{\{X^{\theta}_{t-s_2}\left(\varepsilon,\emptyset,x_{1}\right)\geqslant{}x+y\}}\right],
\end{align*}
and note that
\begin{align*} \left|\widetilde{\Psi}_{s_{2},t}(x_{1})-{\Psi}_{s_{2},t}(x_{1})\right|
&=\left|\mbE\left[\left(\bar{\nu}\left(x_{1},J_{2}\right)-\bar{\nu}\left(X^{\theta}_{t-s_2}\left(\varepsilon,\emptyset,x_{1}\right),J_{2}\right)\right)\mathbf{1}_{\{X^{\theta}_{t-s_2}\left(\varepsilon,\emptyset,x_{1}\right)\geqslant{}x+y\}}\right]\right|\\ &\leqslant\mbE\left[\left|\bar{\nu}\left(x_{1},J_{2}\right)-\bar{\nu}\left(X^{\theta}_{t-s_2}\left(\varepsilon,\emptyset,x_{1}\right),J_{2}\right)\right|\right]\\
&\leqslant\sup_{x,r}|\partial_{1}\bar\nu(x,r)|\mbE\left[\left|X^{\theta}_{t-s_2}(\varepsilon,\emptyset,x_1)-x_1\right|\right].
\end{align*}
Using the facts that $X^{\theta}_{t-s_2}\left(\varepsilon,\emptyset,x_{1}\right)-x_{1}\stackrel{\mathcal{D}}{=}X^{\theta}_{t-s_2}\left(\varepsilon,\emptyset,0\right)$ and $\sup_{x,r}|\partial_{1}\bar\nu(x,r)|<\infty$, the last term above converges to $0$ as $s_2<t\to0$. Therefore, it suffices to study the asymptotic behavior of the expression
\begin{align*}
	\widetilde{T}_{3,2}:&=\mbE\left[\bar\nu\left(\tilde{X}_{s_1}(\varepsilon,\emptyset,x),J_{1}\right)
\left.\mbE\left[\widetilde{\Psi}_{s_{2},t}\left(\hat{X}_{s_2-s_1}(\varepsilon,\emptyset,x_{2})\right)\right]\right|_{x_{2}=\tilde{X}_{s_1}(\varepsilon,\emptyset,x)+\gamma(\tilde{X}_{s_1}(\varepsilon,\emptyset,x),J_1)}\right]\\	&=\mbE\bigg[\bar\nu\left(\tilde{X}_{s_1}(\varepsilon,\emptyset,x),J_{1}\right)\\
&\qquad\quad\times\mbE\left[\left(1-\bar{\nu}\left(X^{\theta}_{t-s_1}\left(\varepsilon,\emptyset,x_{2}\right),J_{2}\right)\right)\mathbf{1}_{\{X^{\theta}_{t-s_1}\left(\varepsilon,\emptyset,x_{2}\right)\geqslant{}x+y\}}\right]\Big|_{x_{2}=\tilde{X}_{s_1}(\varepsilon,\emptyset,x)+\gamma(\tilde{X}_{s_1}(\varepsilon,\emptyset,x),J_1)}\bigg].
\end{align*}
Using again that $\sup_{x,r}|\bar\nu(x,r)|<\infty$, $\sup_{x,r}|\partial_{1}\bar\nu(x,r)|<\infty$ and $\mbE\left[\left|X^{\theta}_{s}\left(\varepsilon,\emptyset,x\right)-x\right|\right]=
\mbE\left[\left|X^{\theta}_{s}\left(\varepsilon,\emptyset,0\right)\right|\right]\to{}0$, uniformly in $x$ as $s\to{}0$,
\begin{align}\label{RDT32}
	\widetilde{T}_{3,2}&=\mbE\left[
	\widetilde{H}\left(\tilde{X}_{s_1}(\varepsilon,\emptyset,x),t-s_{1},x+y\right)
	\right]+o(1),
\end{align}
where
\[
	\widetilde{H}(z,t,q):=\mbE\left[\bar\nu\left(z,J_{1}\right)\left.\left\{\left(1-\bar{\nu}\left(x_{2},J_{2}\right)\right)\mbP\left[X^{\theta}_{t}\left(\varepsilon,\emptyset,x_{2}\right)\geqslant{}q\right]\right\}\right|_{x_{2}=z+\gamma(z,J_1)}\right].
\]
Using a similar argument as in the proof of Lemma \ref{counterpartA.1},
\begin{equation}\label{NEHN}
	\widetilde{H}(z,t,q)=
	\mbE\left[\bar\nu\left(z,J_{1}\right)
	\left(1-\bar{\nu}\left(z+\gamma(z,J_1),J_{2}\right)\right)
	{\bf 1}_{z+\gamma(z,J_1)\geqslant{}q}
	\right]+t\widetilde{\mathcal{R}}^{1}_{t}(z;q),
\end{equation}
where $\sup_{z,q}|\widetilde{\mathcal{R}}^{1}_{t}(z;q)|<\infty$ for $t$ small enough. Furthermore, the first term on the right-hand side of (\ref{NEHN}) can be expressed as
\[
	\int\int_{q}^{\infty}\bar{\nu}(z,\tilde{\gamma}^{-1}(z,r))(1-\bar{\nu}(r,r_2))\widetilde{\Gamma}(r;z)dr\breve{h}_{\varepsilon}(r_2)dr_2,
\]
which is $C^{2}_{b}$ in $z$. The previous fact together with Dynkin's formula as well as (\ref{RDT32})-(\ref{NEHN}) imply that
\begin{align*}
	\widetilde{T}_{3,2}&=\mbE\left[
	\left.\int\int_{x+y}^{\infty}\bar{\nu}(z,\tilde{\gamma}^{-1}(z,r))(1-\bar{\nu}(r,r_2))\widetilde{\Gamma}(r;z)dr\breve{h}_{\varepsilon}(r_2)dr_2\right|_{z=\tilde{X}_{s_1}(\varepsilon,\emptyset,x)}
	\right]+o(1)\\
&=\int\int_{x+y}^{\infty}\bar{\nu}(x,\tilde{\gamma}^{-1}(x,r))(1-\bar{\nu}(r,r_2))\widetilde{\Gamma}(r;x)dr\breve{h}_{\varepsilon}(r_2)dr_2+o(1).
\end{align*}

Using the change of variable $r_1=\tilde{\gamma}^{-1}(x,r)$ and the representation
\begin{align}\label{widetildeGamma}
  \widetilde{\Gamma}(y;z)
  =\partial_y\int\mathbf{1}_{z+\gamma(z,r)<y}\breve{h}_{\varepsilon}(r)dr
  =\frac{h_{\varepsilon}\left(\tilde{\gamma}^{-1}(z,y)\right)}{\lambda_{\varepsilon}}\partial_2\tilde{\gamma}^{-1}(z,y),
\end{align}
together with (\ref{T3inner2})-(\ref{BndT31}), we have
\begin{align}\nonumber
  T_{3}&=\frac{1}{\lambda_{\varepsilon}^{2}}\int\int_{\tilde{\gamma}^{-1}(x,x+y)}^{\infty}\bar{\nu}(x,r_1)(1-\bar{\nu}(\tilde\gamma(x,r_1),r_2))h_{\varepsilon}(r_1)dr_1{h}_{\varepsilon}(r_2)dr_2+o(1)\\
  \label{T3final}&=\frac{1}{\lambda_{\varepsilon}}\int_{\gamma^{-1}(x,y)}^{\infty}\nu_{\varepsilon}(x,r_1)dr_1
  -\frac{1}{\lambda^2_{\varepsilon}}\int_{\gamma^{-1}(x,y)}^{\infty}\nu_{\varepsilon}(x,r_1)\int{\nu}_{\varepsilon}(\tilde\gamma(x,r_1),r_2)dr_2dr_1+o(1).
\end{align}
We next consider the second term, $T_4(x,y)$, in (\ref{2bigjumps}). By Lemma \ref{counterpartA.1},
\begin{align*}
  T_4(x,y)=\mbE\left[H_0\left(\hat{X}_{s_2}(\varepsilon,\{s_1\},x);x+y\right)\right]+O(t).
\end{align*}
Conditioning on $\mathcal{F}_{s_1^-}$,
\begin{align}
  \nonumber{}T_4(x,y)&=\mbE\left[\left.\mbE\left[\mbE\left[H_0\left(\hat{X}_{s_2-s_1}(\varepsilon,\emptyset,v);x+y\right)\right]_{v=z+\gamma(z,J_1)\theta(z,J_1,U_1)}\right]\right|_{z=\tilde{X}_{s_1}(\varepsilon,\emptyset,x)}\right]+O(t)\\
  \nonumber&=\mbE\left[\left.\mbE\left[(1-\bar{\nu}(z,J_1))\mbE\left[H_0\left(\hat{X}_{s_2-s_1}(\varepsilon,\emptyset,z);x+y\right)\right]\right]\right|_{z=\tilde{X}_{s_1}(\varepsilon,\emptyset,x)}\right]\\
  \nonumber&\quad+\mbE\left[\left.\mbE\left[\bar{\nu}(z,J_1)\mbE\left[H_0\left(\hat{X}_{s_2-s_1}(\varepsilon,\emptyset,v);x+y\right)\right]_{v=z+\gamma(z,J_1)}\right]\right|_{z=\tilde{X}_{s_1}(\varepsilon,\emptyset,x)}\right]+O(t)\\
  \label{T41}&=\mbE\left[\left.\mbE\left[(1-\bar{\nu}(z,J_1))\right]\mbE\left[H_0\left(\hat{X}_{s_2-s_1}(\varepsilon,\emptyset,z);x+y\right)\right]\right|_{z=\tilde{X}_{s_1}(\varepsilon,\emptyset,x)}\right]\\
  \label{T42}&\quad+\mbE\left[\left.\int\bar{\nu}(z,\tilde{\gamma}^{-1}(z,r))\mbE\left[H_0\left(\hat{X}_{s_2-s_1}(\varepsilon,\emptyset,r);x+y\right)\right]\widetilde{\Gamma}(r;z)dr\right|_{z=\tilde{X}_{s_1}(\varepsilon,\emptyset,x)}\right]+O(t)
\end{align}
Let us denote the expressions in (\ref{T41}) and (\ref{T42}) by $T_{4,1}$ and $T_{4,2}$, respectively. Next, since $H_{0}(z;q)$ is $C_{b}^{2}$,
\begin{align}\label{asymH0}
  \left|\mbE\left[H_0\left(\hat{X}_{s_2-s_1}(\varepsilon,\emptyset,z);q\right)\right]-H_0(z;q)\right|
  \leqslant\sup_{z,q}|\partial_zH_0(z;q)|\mbE\left[\left|\hat{X}_{s_2-s_1}(\varepsilon,\emptyset,z)-z\right|\right].
\end{align}
Using the facts that $\sup_{z,q}|\partial_zH_0(z;q)|<\infty$ and $\hat{X}_{s_2-s_1}\left(\varepsilon,\emptyset,z\right)-z\stackrel{\mathcal{D}}{=}\hat{X}_{s_2-s_1}\left(\varepsilon,\emptyset,0\right)$, the last term above converges to $0$ uniformly in $z$ as $s_2-s_1<t\to0$. Therefore, recalling the definition of $H_{0}(z;q)=\int^{\infty}_{\gamma^{-1}(z,q-z)}\nu_{\varepsilon}(z,r)dr/\lambda_{\varepsilon}$, we have
\begin{align}
  \nonumber{}T_{4,1}&=\mbE\left[ \left.\mbE\left[(1-\bar{\nu}(z,J_1))\right]\right|_{z=\tilde{X}_{s_1}(\varepsilon,\emptyset,x)}H_0\left(\tilde{X}_{s_1}(\varepsilon,\emptyset,x);x+y\right)\right]+o(1)\\
  \label{T41from}&=\frac{1}{\lambda_{\varepsilon}}\mbE\left[\left.\int_{\gamma^{-1}(z,x+y-z)}^{\infty}\nu_{\varepsilon}(z,r)dr
  \int\left(1-\bar\nu(z,r)\right)\breve{h}_{\varepsilon}(r)dr\right|_{z=\tilde{X}_{s_1}(\varepsilon,\emptyset,x)}\right]+o(1).
\end{align}
Since the expression inside the expectation above has partial derivative with respect to $z$, uniformly bounded in $z$, using again that $\mbE\left[\left|\tilde{X}_{s_1}\left(\varepsilon,\emptyset,x\right)-x\right|\right]=
\mbE\left[\left|\tilde{X}_{s_1}\left(\varepsilon,\emptyset,0\right)\right|\right]\to{}0$ uniformly in $x$ as $s_1<t\to{}0$,
\begin{align}
  \nonumber{}T_{4,1}&=\frac{1}{\lambda_{\varepsilon}}\int_{\gamma^{-1}(x,y)}\nu_{\varepsilon}(x,r)dr
  \int\left(1-\bar\nu(x,r)\right)\breve{h}_{\varepsilon}(r)dr+o(1)\\
  \label{T41final}&=\frac{1}{\lambda_{\varepsilon}}\int_{\gamma^{-1}(x,y)}\nu_{\varepsilon}(x,r)dr
  -\frac{1}{\lambda^2_{\varepsilon}}\int_{\gamma^{-1}(x,y)}\nu_{\varepsilon}(x,r)dr\int\nu_{\varepsilon}(x,r)dr+o(1).
\end{align}

We next consider the term $T_{4,2}(t)$ in (\ref{T42}). Using (\ref{asymH0}) again and an argument similar to that from (\ref{T41from}) to (\ref{T41final}), we deduce as follows.
\begin{align}
  T_{4,2}&=\mbE\left[\left.\int\bar{\nu}(z,\tilde{\gamma}^{-1}(z,r))H_0(r;x+y)\widetilde{\Gamma}(r;z)dr\right|_{z=\tilde{X}_{s_1}(\varepsilon,\emptyset,x)}\right]+o(1)\nonumber \\
  &=\frac{1}{\lambda_{\varepsilon}}\int\bar{\nu}(x,\tilde{\gamma}^{-1}(x,r))\int^{\infty}_{\gamma^{-1}(r,x+y-r)}\nu_{\varepsilon}(x,r_2)dr_2\widetilde{\Gamma}(r;x)dr+o(1)\nonumber \\
  &=\frac{1}{\lambda^2_{\varepsilon}}\int\nu_{\varepsilon}(x,r_1)\int^{\infty}_{\gamma^{-1}(x+\gamma(x,r_1),y-\gamma(x,r_1))}\nu_{\varepsilon}(x,r_2)dr_2dr_1+o(1),\label{T42final}
\end{align}
where, in the last equality, we use the change of variable $r_1=\tilde{\gamma}^{-1}(x,r)$ and the representation (\ref{widetildeGamma}) of $\widetilde{\Gamma}(y;z)$. Summing up (\ref{T41final}) and (\ref{T42final}), we have
\begin{equation}\label{T4final}
\begin{aligned}
  T_4(x,y)&=\frac{1}{\lambda_{\varepsilon}}\int^{\infty}_{\gamma^{-1}(x,y)}\nu_{\varepsilon}(x,r)dr
  -\frac{1}{\lambda^2_{\varepsilon}}\int^{\infty}_{\gamma^{-1}(x,y)}\nu_{\varepsilon}(x,r)dr\int\nu_{\varepsilon}(x,r)dr\\
  &\quad+\frac{1}{\lambda^2_{\varepsilon}}\int\nu_{\varepsilon}(x,r_1)\int^{\infty}_{\gamma^{-1}(x+\gamma(x,r_1),y-\gamma(x,r_1))}\nu_{\varepsilon}(x,r_2)dr_2dr_1+o(1).
\end{aligned}
\end{equation}
Therefore, summing up $T_3(x,y)$ and $T_4(x,y)$ from (\ref{T3final}) and (\ref{T4final}), $\mbP\left[X^{\theta}_t(\varepsilon,\{s_1,s_2\},x)\geqslant{}x+y\right]=H_{2,0}(x;y)+o(1)$ with
\begin{equation}\label{2JExp}
\begin{aligned}
  H_{2,0}(x;y)
  &=\frac{2}{\lambda_{\varepsilon}}\int^{\infty}_{\gamma^{-1}(x,y)}\nu_{\varepsilon}(x,r)dr
  -\frac{1}{\lambda^2_{\varepsilon}}\int_{\gamma^{-1}(x,y)}^{\infty}\nu_{\varepsilon}(x,r_1)\int\nu_{\varepsilon}(\tilde\gamma(x,r_1),r_2)dr_2dr_1\\
  &\quad-\frac{1}{\lambda^2_{\varepsilon}}\int^{\infty}_{\gamma^{-1}(x,y)}\nu_{\varepsilon}(x,r)dr\int\nu_{\varepsilon}(x,r)dr\\
  &\quad+\frac{1}{\lambda^2_{\varepsilon}}\int\nu_{\varepsilon}(x,r_1)\int^{\infty}_{\gamma^{-1}(x+\gamma(x,r_1),y-\gamma(x,r_1))}\nu_{\varepsilon}(x,r_2)dr_2dr_1,
\end{aligned}
\end{equation}
and, due to (\ref{N2int}),
\begin{align}\label{ConN2}
\mbP\left[\left.X^{\theta}_{t}(\varepsilon,x)\geqslant{}x+y\right|N_{t}^{\varepsilon}=2\right]
=H_{2,0}(x;y)+o(1).
\end{align}

\noindent{}Substituting (\ref{ConN1}) and (\ref{ConN2}) into (\ref{conditional}) and using the fact that $e^{-\lambda_{\varepsilon}t}=1-\lambda_{\varepsilon}t+\frac{(\lambda_{\varepsilon}t)^2}{2}+O(t^3)$ as $t\downarrow0$, we obtain
\begin{equation*}
  \mbP[X_t\geqslant{}x+y]=tP_1(x,y)+\frac{t^2}{2}P_2(x,y)+O(t^3),
\end{equation*}
where
\begin{align*}
  P_1(x,y)&=\lambda_{\varepsilon}H_{0,0}(x;y)
  =\int_{\gamma^{-1}(x,y)}^{\infty}\nu_{\varepsilon}(x,r)dr,\\
  P_2(x,y)&=-2\lambda^2_{\varepsilon}H_{0,0}(x;y)
  +\lambda_{\varepsilon}[H_{0,1}(x;y)+H_{1,0}(x;y)]
  +\lambda^2_{\varepsilon}H_{2,0}(x;y).
\end{align*}
Some further simplifications lead to the expressions of $P_{1}(x,y)$ and $P_{2}(x,y)$ stated in the statement of the theorem.
\hfill$\blacksquare$

\subsection{Corollary \ref{Ppoundexpansion}}\label{ProofPpoundexpansion}

\textbf{Proof: }
Let us first note that, due to the Condition (\ref{conditiongamma}-\ref{gamma0=0}),
\begin{align*}
	|\gamma(x,r)|\leqslant{}cr,
\end{align*}
for all $x$ and $r$, where $c$ is defined as in Condition (\ref{expgamma}).
Next, define $g^{\#}(r):=g(r)e^{cr}$ and note that, in view of the Condition (\ref{expgamma}), $h^{\#}(r):=g^{\#}(r)|r|^{-\alpha-1}$ is a valid L\'evy density. Furthermore, $g^{\#}$ clearly satisfies the other requirements of Condition (\ref{conditionh}), while
\[
	\bar\nu^{\#}(x,r):=\frac{\nu^{\#}(x,r)}{h^{\#}(r)}=\frac{e^{\gamma(x,r)}\nu(x,r)}{e^{cr}g(r)|r|^{-\alpha-1}}=e^{\gamma(x,r)-cr}\bar\nu(x,r),
\]
can be readily seen to meet all the requirements of Condition (\ref{conditionnu}). The result is then a consequence of Theorem \ref{tailexpansion}.
\hfill$\blacksquare$

\subsection{Lemma \ref{smallJ}}\label{pfsmallJ}
\begin{proof}
From the infinitesimal generators (\ref{smallJGenerator})-(\ref{InfGenTildeL}), we can identify the semimartingale characteristics  $(B,C,\vartheta)$ and $(B^{\varepsilon},C^{\varepsilon},\vartheta^{\varepsilon})$ of  $\{{X}_t\}_{t\geqslant0}$ and $\{\widetilde{X}_t^{\varepsilon}\}_{t\geqslant0}$, respectively,
relative to the truncation function $r{\bf 1}_{\{|r|\leqslant{}1\}}$:
\begin{align}
\label{SemiCharacterX} B_t=\int_0^t\hat{b}\left({X}_s\right)ds,
\quad C_t&=\int_0^t{}\sigma^2\left(X_s\right)ds,
\quad \vartheta(dt,dr)=K\left(X_{t^-},dr\right)dt,\\
\label{SemiCharacterEp} B_t^{\varepsilon}=\int_0^t\hat{b}\left(\widetilde{X}^{\varepsilon}_s\right)ds,
\quad C_t^{\varepsilon}&=\int_0^t{}\hat\sigma_{\varepsilon}^2\left(\widetilde{X}^{\varepsilon}_s\right)ds,
\quad \vartheta^{\varepsilon}(dt,dr)=\widetilde{K}_{\varepsilon}\left(\widetilde{X}^{\varepsilon}_{t^-},dr\right)dt.
\end{align}
where $K\left(x,A\right):=\int {\bf 1}_{A}\left(\gamma\left(x,r\right)\right)\nu\left(x,r\right)dr$, $\bar{b}_{\varepsilon}(x)$. We also define $\bar{b}_{\varepsilon}(x)$ and $\bar{\sigma}_{\varepsilon}^2(x)$ as 
\begin{align}\label{tildeb}
	\bar{b}_{\varepsilon}(x)&:=\hat{b}(x)+\int{}{\widetilde{K}_{\varepsilon}}(x,dr)r{\bf 1}_{\{|r|>{}1\}}=\hat{b}(x)+\int{}\gamma(x,r){\bf 1}_{\{|\gamma(x,r)|>{}1\}}\nu(x,r){\bf 1}_{\{|\gamma(x,r)|>\varepsilon\}}dr\\
	\label{tildesigma}
	\bar{\sigma}_{\varepsilon}^2(x)
	&:=\tilde\sigma_{\varepsilon}^2(x)+\int{}{\widetilde{K}_{\varepsilon}}(x,dr)r^2=\sigma^{2}(x)+\hat\sigma_{\varepsilon}^2(x)+\int{}\gamma^{2}(x,r)\nu(x,r){\bf 1}_{\{|\gamma(x,r)|>\varepsilon\}}dr.
\end{align}
Clearly, for $\varepsilon\in(0,1)$, $\bar{b}_{\varepsilon}=\hat{b}(x)+\int{}K(x,dr)r{\bf 1}_{\{|r|>{}1\}}$ and
$\bar\sigma_{\varepsilon}^{2}(x)=\sigma^{2}(x)+\int{}K(x,dr)r^2$,
in light of the definition of $\hat{\sigma}_{\varepsilon}^{2}$. Also,
\[
	\int {\widetilde{K}_{\varepsilon}}(x,dr)g(r)=\int{}g(\gamma(x,r))\nu(x,r){\bf 1}_{\{|\gamma(x,r)|>\varepsilon\}}dr\stackrel{\varepsilon\to{}0}{\longrightarrow}\int{}g(\gamma(x,r))\nu(x,r)dr=\int K(x,dr)g(r),
\]
for any bounded continuous function $g$ that vanishes in a neighborhood of the origin. Therefore, by \cite[Theorem IX.4.15]{Jacod}, $\{\widetilde{X}_t^{\varepsilon}\}_{t\geqslant0}\xrightarrow{\mathcal{D}}\{X_t\}_{t\geqslant0}$ provided that the uniqueness hypothesis \cite[IX.4.3]{Jacod} holds true.
By \cite[Theorem III.2.26]{Jacod}, the uniqueness requirement stated in \cite[IX.4.3]{Jacod} is equivalent to the weak uniqueness of the solution for the SDE defining $X$. This fact was established in Lemma \ref{Y0diffeomorphism}. Therefore, we conclude the convergence result claimed in the lemma by \cite[Theorem IX.4.15]{Jacod}.
\end{proof}

\sectionfont{}

\end{document}